\documentclass[11pt,reqno]{amsart}
\usepackage{amsbsy,amsfonts,amsmath,amssymb,amscd,amsthm}
\usepackage{mathrsfs}
\usepackage{subfigure,enumitem,color,graphicx}
\usepackage[a4paper,text={6.1in,8in},centering,includefoot,foot=1in]{geometry}
\usepackage{pxfonts}
\usepackage[colorlinks=true,linkcolor=blue,citecolor=green]{hyperref}
\usepackage{srcltx}

\small\normalsize
\newtheorem{theo}{Theorem}[section]
\newtheorem{prop}[theo]{Proposition}
\newtheorem{lemma}[theo]{Lemma}

\newtheorem{mainthm}{Theorem}

\newtheorem{claim}{Claim}
\newtheorem{addendum}[theo]{Addendum}
\newtheorem{defi}[theo]{Definition}

\newtheorem{rem}[theo]{Remark}

\newcommand{\eqdef}{\stackrel{\scriptscriptstyle\rm def}{=}}
\DeclareMathOperator{\IFS}{IFS} 

\begin{document}
\title[IFS on the circle close to rotations]{Dynamics of iterated function systems on the circle close to rotations}

\author[Barrientos]{Pablo G. Barrientos}
\address{\footnotesize \centerline{Departamento de Matem\'atica PUC-Rio},
\centerline{Marqu\^es de S\~ao Vicente 225, \\
G\'avea, Rio de Janeiro 225453-900, Brazil}}
\email{barrientos@mat.puc-rio.br}

\author[Raibekas]{Artem Raibekas}
\address{\footnotesize \centerline{Instituto de Matem\'atica e Estat\'istica, UFF}
\centerline{Rua M\'ario Santos Braga s/n - Campus Valonguinhos,
Niter\'oi, Brazil}} \email{artem@mat.uff.br}

\maketitle

\begin{abstract}
We study the dynamics of iterated function systems generated by a
pair of circle diffeomorphisms close to rotations in the
$C^{1+\mathrm{bv}}$-topology. We characterize the obstruction to
minimality and describe the limit set. In particular, there are no
invariant minimal Cantor sets, which can be seen as a
Denjoy/Duminy type theorem for iterated systems on the
circle.
\end{abstract}

\begin{flushright}
\textit{D\'edi\'e \`a G. Duminy}
\end{flushright}

\section{Introduction}

The rotation number is a classical tool in the study of dynamics of a single
diffeomorphism of the circle. When the rotation number is rational, the map
has periodic points and if it is irrational then either each orbit is dense
on the circle or there is an invariant (minimal) Cantor set. Consider now a group of orientation-preserving homeomorphisms of the circle.
The above classification
can be extended to the group case where the orbit of a point is the group action on this point. Then, there can occur only one of the following
three options~\cite{Gh01, Na11}:
existence of a finite orbit, every orbit is dense on the circle, or there exists a unique minimal
Cantor set invariant by the group action. Minimal here means that all the points in the invariant set have dense orbits.

For iteration of single diffeomorphisms, the well-known theorem by Denjoy~\cite{Denjoy} implies  that in the $C^2$-topology there cannot exist minimal invariant Cantor sets,
while there are counter-examples in the $C^1$ class.
Motivated by the study of co-dimension one foliations, Duminy proved a
Denjoy-like theorem for group actions~\cite{Du}. Under the extra assumptions that the group generators are close to rotations
in the $C^{1+\mathrm{bv}}$-topology (i.e. the $C^1$-class with bounded variation derivatives) and \emph{finiteness of periodic points} for one of the generators,
the minimal invariant Cantor set does not exist. See~\cite{Na11,Na04} for details of the proof.

In this paper, motivated by the dynamics of partially hyperbolic skew-products~\cite{GI99,BDV05,NP11,BKR}, we
study actions of finitely generated semigroups of diffeomorphisms on the circle, which can be viewed as iterated function systems.
We extend the theorems of Denjoy and Duminy to the semigroup case and describe the limit set of possible orbits.

In order to state the main results, first a couple of definitions.

An \emph{iterated function system} (IFS),
 generated by a finite family of diffeomorphisms
 $\Phi=\{\phi_1,\ldots,\phi_k\}$ of the circle $\mathbb{S}^1$,
 is the set $\IFS(\Phi)$ of all nonempty
possible finite compositions of diffeomorphisms $\phi_i \in \Phi$.
That is, the semigroup
generated by the compositions of
$\phi_1,\ldots,\phi_k$.

The \emph{orbit} of $x$ for $\IFS(\Phi)$ is the action of IFS
over the point $x$, i.e.,
$$
\mathrm{Orb}_\Phi^+(x)\eqdef \{h(x): \, h \in \IFS(\Phi)\} \subset
\mathbb{S}^1.
$$
Denote by $\mathrm{Per}(\IFS(\Phi))$ the \emph{set of periodic
points} of $\IFS(\Phi)$, which is the set of $x\in \mathbb{S}^1$
such that $h(x)=x$ for some
$h \in \mathrm{IFS}(\Phi)$. 
We say the circle $\mathbb{S}^1$ is minimal for $\IFS(\Phi)$ if
the orbit of every point is dense.

Let $f$ be a $C^{1+\mathrm{bv}}$-diffeomorphism of the circle.
From the Denjoy Theorem the following statements are equivalent:
$\mathbb{S}^1$ is minimal for $\IFS(f)$ and $f$ does not have
periodic points. When the number of generators of the IFS is at
least 2, the periodic points are no longer the unique obstruction
to minimality. Indeed, this role is now played by the
$ss$-intervals which are compact intervals whose endpoints are
consecutive attracting fixed points of different generators of
$\mathrm{IFS}(\Phi)$.
See Definition~\ref{d:**intervals} for a formal definition. \\[-0.5cm]

\begin{mainthm}[Obstruction to minimality]
\label{thmA}
There exists $\varepsilon>0{.}38$ such that if $f_0$ and $f_1$ are diffeomorphisms of the circle
with periodic points of period $n_0$ and $n_1$ respectively, $\varepsilon$-close to the rotations in the $C^{1+\mathrm{bv}}$-topology and with
no periodic points in common, then the following conditions are equivalent:
\begin{itemize}
\item $\mathbb{S}^1$ is minimal for $\IFS(f_0^{n_0},f_1^{n_1})$,
\item $\mathbb{S}^1$ is minimal for $\IFS(f_0^{-n_0},f_1^{-n_1})$,
\item there are no $ss$-intervals for $\IFS(f_0^{n_0},f_1^{n_1})$.
\end{itemize}
    When this occurs, both the hyperbolic attracting and the hyperbolic repelling periodic points of $\mathrm{IFS}(f_0^{n_0},f_1^{n_1})$ are dense in $\mathbb{S}^1$. \\[-0.5cm]
\end{mainthm}

Under the same assumptions of Theorem~\ref{thmA}, the first
theorem of Duminy~\cite{Na11, Na04} can be restated as,
$\mathbb{S}^1$ is minimal for the group action and the hyperbolic
periodic points are dense. See Theorem~\ref{t.Duminy} in this
article.

If we assume the generators have only hyperbolic periodic points,
then Theorem~\ref{thmA} is robust in the following sense. There
exists an open set $\mathcal{U}\subset
\mathrm{Diff}^{1+\mathrm{bv}}(\mathbb{S}^1)\times
\mathrm{Diff}^{1+\mathrm{bv}}(\mathbb{S}^1)$ with $(f_0,f_1)\in
\mathcal{U}$ such that for every pair $(g_0, g_1) \in \mathcal{U}$
it holds that $\mathbb{S}^1$ is minimal for
$\IFS(g_0^{n_0},g_1^{n_1})$ if and only if there are no
$ss$-intervals for $\IFS(f_0^{n_0},f_1^{n_1})$. Note that examples
of robustly minimal IFS on the circle
were given in~\cite{GI99,GI00,KN04}.

Actually, the above robustness of Theorem~\ref{thmA} is valid for
perturbations of $(f_0,f_1)$ in the $C^1$-topology, i.e., in
$\mathrm{Diff}^{1}(\mathbb{S}^1)\times
\mathrm{Diff}^{1}(\mathbb{S}^1)$. This follows from the work of
\cite{BFS} using the fact that the semigroup action is expanding
(see Remark~\ref{exp}). The action of a semigroup $\Gamma$ of
orientation-preserving diffeomorphisms of the circle is said to be
\emph{expanding} if for any point $x\in \mathbb{S}^1$ there exists
$g\in \Gamma$ with $Dg^{-1}(x)>1$. Note that in the same
paper~\cite{BFS} is shown the robust ergodicity of expanding
semigroup actions (in the $C^{1+\alpha}$ topology), where the
proof is based on a similar result for group actions of the
circle~\cite{Na04}.

The existence of $ss$-intervals and the non-minimality of the circle leads us to attempt in
understanding the indecomposable pieces of global dynamics. Thus, we introduce the notion of a limit set for the orbits of an IFS.

The \emph{forward} or \emph{$\omega$-limit of $x\in \mathbb{S}^1$
for  $\IFS(\Phi)$} is the set
$$
\omega_\Phi(x)\eqdef
\{y: \, \exists \, (h_n)_n \subset \IFS(\Phi) 
\ \text{such that}\, \displaystyle\lim_{n\to \infty} h_n\circ\dots\circ h_1(x)=y \},
$$
while the \emph{$\omega$-limit of $\IFS(\Phi)$} is
$$
\omega(\IFS(\Phi)) \eqdef \mathrm{cl}\big(\{ y: \ \exists \, x \in M \text{ such that } \ y \in \omega_\Phi(x)\}\big),
$$
where ''$\mathrm{cl}$´´ denotes the closure of a set. Similarly,
the \emph{backward} or \emph{$\alpha$-limit of $\IFS(\Phi)$} is
defined as $
  \alpha(\IFS(\Phi))\eqdef \omega(\IFS(\Phi^{-1}))
$
where $\Phi^{-1}=\{\phi_1^{-1},\ldots,\phi_k^{-1}\}$.
From the backward and forward limit, we define the \emph{limit set of $\IFS(\Phi)$} as
$$
   L(\IFS(\Phi))\eqdef
   \omega(\IFS(\Phi)) \cup  \alpha(\IFS(\Phi)).
$$

The following result, in particular, shows the spectral
decomposition of the limit set of IFS generated by Morse-Smale
diffeomorphisms of the circle. That is, the limit set can be
written as a finite union of disjoint (closed) maximal transitive
sets. A (not necessarily IFS invariant) set $\Lambda \subset
\mathbb{S}^1$ is \emph{transitive} for $\IFS(\Phi)$ if
$$
\Lambda \subset \overline{\mathrm{Orb}^+_\Phi(x)} \quad  \text{for some $x\in\Lambda$.}
$$
The transitive set is \emph{maximal} if there is no other transitive set $\Omega$ with $\Lambda\subset\Omega$.

\begin{mainthm}[Spectral decomposition]
\label{thmB} There exists $\varepsilon>0{.}30$ such that, under
the assumptions of Theorem~\ref{thmA} and if \, $\mathbb{S}^1$ is
not minimal for $\IFS(f_0^{n_0},f_1^{n_1})$, then
$$
   L(\IFS(f_0^{n_0},f_1^{n_1}))=\bigcup K_i \quad \text{and} \quad
   K_i\subset \overline{\mathrm{Per}(\IFS(f_0^{n_0},f_1^{n_1}))}
$$
where each $K_i$ is a 
compact, maximal transitive set for $\IFS(f_0^{n_0},f_1^{n_1})$.
Moreover, each compact set $K_i$ is either \begin{itemize}
\item a single fixed point of $f_0^{n_0}$ or $f_1^{n_1}$, or
\item a $**$-interval  for $\IFS(f_0^{n_0},f^{n_1}_1)$
with $**\in \{ss,su,uu\}$.\footnote{For the precise definition see Definition~\ref{d:**intervals}.}
\end{itemize}
In particular, if $f_0$ and $f_1$ are Morse-Smale diffeomorphisms then the above union is finite and the sets $K_i$ are pairwise disjoints.
%
\end{mainthm}

Observe that the spectral decomposition of Theorem~\ref{thmB} is given for the IFS generated by $f_0^{n_0}$ and $f_1^{n_1}$.
It remains a question to describe the possible limit set of $\IFS(f_0,f_1)$, and if there is a spectral decomposition what are the pieces.
An extension of the problem is to prove a spectral decomposition type result when the generators
are far from rotations. Observe that in this case Cantor sets may appear in the decomposition.

With respect to dynamical decomposition of IFS from the ergodic perspective, let us mention
the works of \cite{KV} and \cite{Viana2012}, where is shown the existence of finitely many
stationary or SRB measures for step skew-products over the interval or the circle. On the contrary, for non-Abelian
group actions of the circle in \cite{DKN07} the uniquenenss of stationary measures is proven.

A set $\Lambda$ is (forward) invariant for $\IFS(\Phi)$ if $\phi_i(\Lambda)\subset \Lambda$ for all $i=1,\ldots,k$. We say that a closed
invariant set $\Lambda$ is \emph{minimal} for $\IFS(\Phi)$ if
every point of $\Lambda$ has dense orbit in $\Lambda$
or equivalently if
$$
\Lambda=\phi_1(\Lambda)\cup \dots  \cup \phi_k(\Lambda)=\overline{\mathrm{Orb}^+_\Phi(x)}  \quad \text{for all $x\in\Lambda$}.
$$

Similarly to the case of group actions, a closed minimal invariant set $\Lambda$ for IFS has to satisfy one of the following:
 $\Lambda$ is a finite orbit, has non-empty interior, or is a Cantor set, see Theorem~\ref{p.IFS-tricotomi}.
Observe that apriori minimal invariant Cantor sets could exist in the pieces of the spectral decomposition of Theorem~\ref{thmB}.
The next result shows that this cannot occur and can be thought of
as a Denjoy-type result for the semigroup action.

\begin{mainthm}[Denjoy for IFS]
\label{thmC} There exists $\varepsilon>0{.30}$ such that under the
assumptions of Theorem~\ref{thmA} and if \, $\mathbb{S}^1$ is not
minimal for $\IFS(f_0,f_1)$, the only minimal invariant closed
sets for $\IFS(f_0^{n_0},f_1^{n_1})$ are the $ss$-intervals.

Moreover, there are no minimal invariant Cantor sets for
$\IFS(f_0,f_1)$.
\end{mainthm}

In the case of iterations of a single diffeomorphism on the circle, the Denjoy examples~\cite[Chapter X]{Her79}
show the existence of minimal invariant Cantor sets in the
$C^1$-topology. Thus one can expect that Theorem~\ref{thmC} above fails in the $C^1$-topology.
The recent works of Shinohara~\cite{S12, S13} give examples of $\IFS(f_0,f_1)$
with $f_0$ and $f_1$  close to the identity in the $C^1$-topology and forming an $ss$-interval,
such that the minimal set is a Cantor set.


The paper is organized as follows.
Firstly, in Section~\ref{s:duminy} we state and give a proof of a relevant result of Duminy. This is necessary
in order to understand the next section, which is a generalization of Duminy's Theorem
using a new notion of cycles for IFS and concludes with the proof of Theorem~\ref{thmA}. In Section~\ref{s:spectral-decomposition}
we show a Spectral Decomposition Theorem on the real line, and afterwards prove Theorem~\ref{thmB}.
The trichotomy on the shape of minimal invariant sets is taken up in the next section and finally
Theorem~\ref{thmC} and Duminy's Theorem under the assumptions of Theorem~\ref{thmA} are shown.

\noindent \textbf{Notation:} In what follows we will use the
following notation, where $f_0$ and $f_1$ are diffeomorphisms of
the circle with periodic points of period $n_0$ and $n_1$
respectively:
\begin{align*}
\Phi&=(f_0,f_1), & \Phi^{-1}&=(f_0^{-1},f_1^{-1}), &  \Phi^n&=(f_0^{n_0},f_1^{n_1}), & 
\Phi^{-{n}}&=(f_0^{-n_0},f_1^{-n_1}).
\end{align*}
For a map $f$ defined on an interval $I$, we write $f>\mathrm{id}$
in $I$ if $f(x)>x$ for all $x\in I$. Also we denote by $|I|$ the
length of this interval.
\medbreak

\section{A result of Duminy}
\label{s:duminy}
In this section we will show an important result of Duminy which can be easily deduced from~\cite{Na04, Na11}.
For the sake of a self-contained paper and to indicate the techniques used to derive the stronger result in Section~\ref{ss:cycle},
the complete proof is given. Before that, the formal definition of $**$-intervals is introduced.

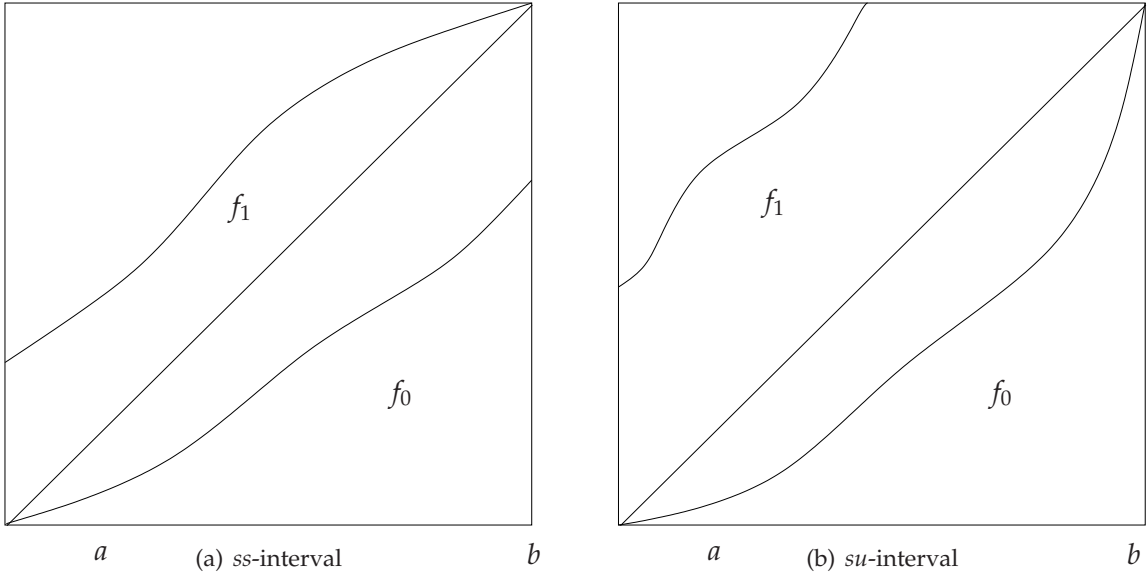
\begin{figure}
\centering
\begin{picture}(425,210)
\subfigure[$ss$-interval]{\label{Kss}%
\includegraphics[height=0.45\textwidth,width=0.45\textwidth]{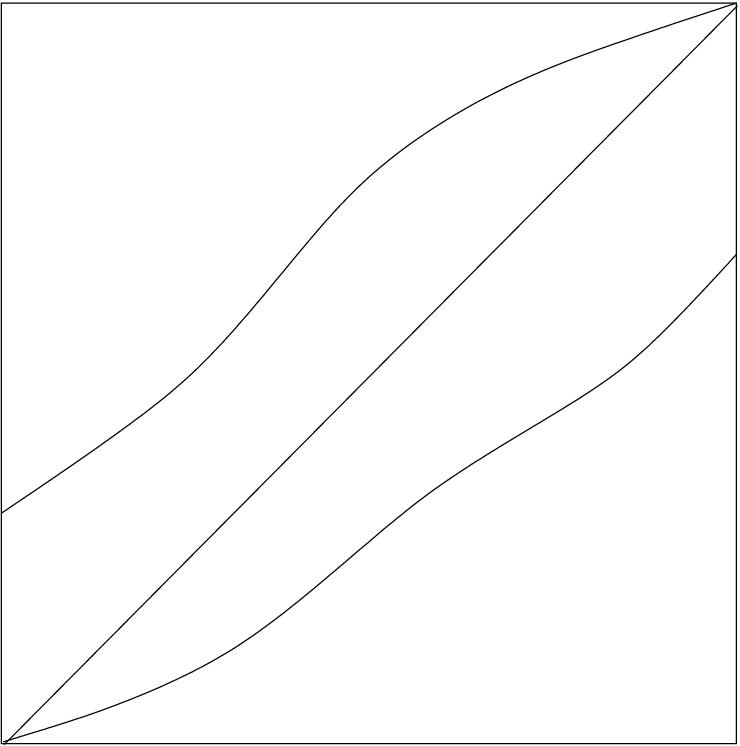}}
\put(-120,110){\makebox(20,20){$f_1$}}
\put(-60,40){\makebox(20,20){$f_0$}}
\put(-172,-20){\makebox(20,20){$a$}}
\put(-10,-20){\makebox(20,20){$b$}} \hspace{1cm}
\subfigure[$su$-interval]{\label{Ksu}%
\includegraphics[height=0.45\textwidth,width=0.45\textwidth]{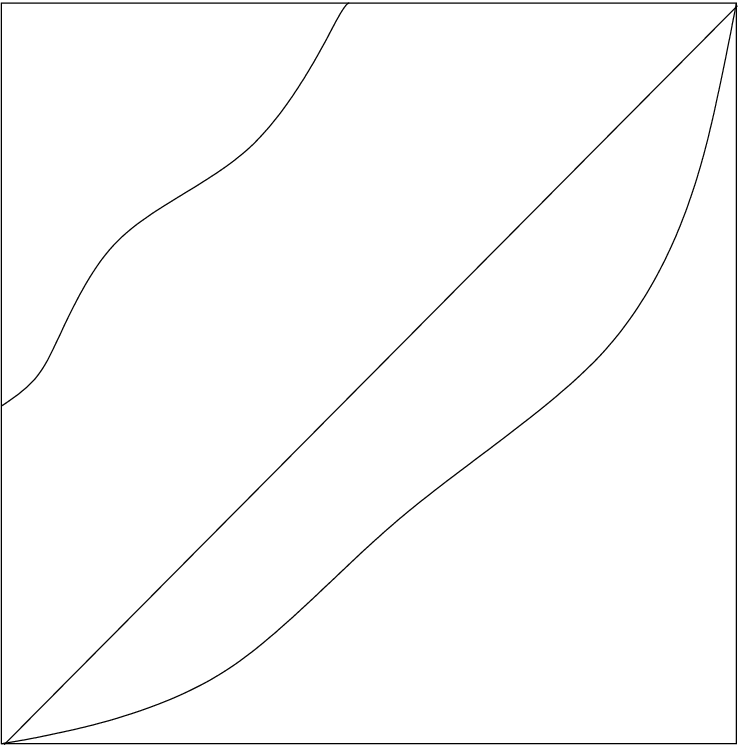}}
\put(-150,112){\makebox(20,20){$f_1$}}
\put(-65,40){\makebox(20,20){$f_0$}}
\put(-172,-20){\makebox(20,20){$a$}}
\put(-15,-20){\makebox(20,20){$b$}}
\end{picture}
\vspace{0.5cm}
\caption{Examples of ${**}$-intervals}
\label{**-intervalos}
\end{figure}

\begin{defi}[$**$-intervals]
\label{d:**intervals}
Consider $f_0$, $f_1$ orientation preserving homeomorphims on the real line. Let $[a,\,b]$ be a compact interval such that
$$
\mathrm{Fix}(f_i)\cap (a,\,b) =
\emptyset \quad \text{and} \quad
   [a,\,b]\subset f_0([a,\,b])\cup f_1([a,\,b]).
$$
We say that $[a,\,b]$ is a \emph{$**$-interval} for $\IFS(\Phi)$  
with $** \in \{
ss, su\}$  when $a$ and $b$ satisfy the additional properties (see Figure~\ref{**-intervalos}):
\begin{itemize}
\item \emph{$ss$-interval}:
$a$ and $b$ are, respectively, attracting fixed points of the
restriction of $f_0$ and $f_1$ to $[a,\,b]$  and $f_0(b)\not=b$,
$f_1(a)\not=a$.
\item \emph{$su$-interval}: $a$ and $b$ are either, an attractor-repeller or  a repeller-attractor pair for the same map restricted to $[a,\,b]$, say an attractor-repeller pair for $f_0$. In this case, we ask that $f_1>\mathrm{id}$ in $[a,\,b]$ and $f_1([a,\,b])\cap (a,\,b)\neq \emptyset$.
\end{itemize}
A \emph{$uu$-interval} for $\IFS(\Phi)$ is a $ss$-interval for
$\IFS(\Phi^{-1})$. An unbounded interval $[a,\,\infty]$ is
\begin{itemize}
\item $s$-interval for $\IFS(\Phi)$:  if $a$ is an attracted fixed point of the restriction of a map to $[a,\,\infty)$, say $f_0$, satisfying
    $f_0<\mathrm{id}$ in $(a,\,\infty)$  and
    $f_1>\mathrm{id}$ in $[a,\,\infty)$.
\end{itemize}
A \emph{$u$-interval} for $\IFS(\Phi)$ is a $s$-interval for $\IFS(\Phi^{-1})$.
Analogously, the unbounded case $[-\infty,\,b]$ is defined in the same manner.  
\end{defi}

Let $I$ be an interval on the real line or on $\mathbb{S}^1$.
Consider an orientation preserving $C^1$-map $f$ of $I$ such that
$Df(x)\not= 0 $ for all $x\in I$. The non-negative number
\begin{equation*}
   \mathrm{Dist}(f,I)= \sup_{x,y\in I} \log \frac{Df(x)}{Df(y)}
\end{equation*}
is called \emph{distortion constant} of $f$ in $I$ or \emph{maximal variation} of $\log Df$ in $I$. We say that $f$ belongs to the $C^{1+\mathrm{bv}}$ class on $I$
if it has a bounded distortion constant in~$I$. Observe that, an application of the Mean Value Theorem shows that $C^2$-maps with
non-zero first derivative are $C^{1+\mathrm{bv}}$.

It is not difficult to see that every orientation preserving $C^{1+\mathrm{bv}}$-map $f$  on the circle satisfies
$$
  e^{-V_f} \leq Df(x) \leq e^{V_f}  \quad \text{and} \quad  |f(x)-(x+\rho(x))|<1-e^{-V_f} \quad \text{for all $x\in \mathbb{S}^1$}
$$
where $\rho(f)$ is the rotation number of $f$ and
$V_f=\mathrm{Dist}(f,\mathbb{S}^1)$ (see~\cite{Na04} for details).
Therefore, if the distortion constant of $f$ is small enough, then
$f$ is close to a rotation. We will say that $f$ is
\emph{$\varepsilon$-close to rotation in the
$C^{1+\mathrm{bv}}$-topology} if $V_f \leq \varepsilon$.

Let $f_0$ and $f_1$ be orientation preserving $C^{1+\mathrm{bv}}$-diffeomorphisms on the real line.
In what follows, we will study $\IFS(\Phi)$ with $f_0$ and $f_1$ restricted to a ${**}$-interval $K^{**}=[a,\,b]$ for  $** \in \{ss, s, su\}$.
In the case of ${uu}$-intervals and $u$-intervals for $\IFS(\Phi)$ the same results hold for $\IFS(\Phi^{-1})$. Observe that
$$
  K^{**}=f_0(K^{**})\cup f_1(K^{**}) \quad \text{for $**\in\{ss,s\}$}
$$
and then, it is an invariant set by the action of $\IFS(\Phi)$. However the above equality does not hold
for a ${su}$-interval.

We will be interested in stating conditions under which
every point in a ${**}$-interval has a dense orbit for $\IFS(\Phi)$.
In the case of a $su$-interval and a $s$-interval
notice that one of the endpoints of $K^{su}$ and $K^{s}$ can
never have a dense orbit.

To unify the notation, we will define
an $**$-interval to be \emph{minimal} for $\IFS(\Phi)$ when
$$
K^{**}\subset \overline{\mathrm{Orb}^+_\Phi(x)} \ \ \ \text{for all $x\in K^{**}$}.
$$

The details of the following result of Duminy were worked out by Navas~\cite{Na04}, who
reinterpreted the main arguments of the work of Duminy in more dynamical terms, introducing
certain return maps and showing they are expanding.

\begin{theo}[Duminy's Lemma]
\label{t.duminy_lemma}
Assume $f_0 < \mathrm{id}$ in $(a,b)$ and  suppose that there exists $\varepsilon>0$ such that
\begin{align*}
|Df_0(x)-1|<\varepsilon \quad \text{for all } x\in (a,\,f_0^{-1}f_1(a)) \quad 
&\text{and} \quad (1-\varepsilon)\varepsilon^{-1} e^{-V}>1
\end{align*}
where $V=(V_0 + V_1)|a-f_0^{-1}f_1(a)|$ being $V_0$ and $V_1$   constants satisfying
\begin{align*}
  \mathrm{Dist}(f_0, I) \leq V_0|I|  \quad \text{and} \quad \mathrm{Dist}(f_1,J) \leq V_1 |J|
\end{align*}
for all fundamental domains $I$ of $f_0$ and $J$ of $f_1$ in $(a,\,f_0^{-1}f_1(a)]$.
Then
\begin{equation*}
K^{**}  \subset
\overline{\mathrm{Per}(\mathrm{IFS}(\Phi))} \quad \text{and} \quad
K^{**} \subset \overline{\mathrm{Orb^+_\Phi}(x)} \ \ \text{for all} \
x\in K^{**}.
\end{equation*}
\end{theo}

Notice that from definition of ${**}$-intervals for $\IFS(\Phi)$ with $**\in \{ss,su,s\}$, the \emph{overlap condition} is verified, that is,
$
  f_0(K^{**})\cap f_1(K^{**})\not = \emptyset.
$
This condition implies that
$$
A=(f_1(a),\,f_0^{-1}f_1(a)] \subset [a,\,b].
$$
Next, we will define a first return map $\mathcal{R}$ over the fundamental domain $A$.

For each $x\in A$ let $m(x)\geq 1$ be the smallest positive number such that $f_1^{-m(x)}(x)\not\in A$ and let $n(x)$ be the first time for which $f_1^{-m(x)}(x)$ returns to $A$ by iterations of $f_0^{-1}$. Then we can define the first return map $\mathcal{R}$ in the following way
\begin{equation*}
\mathcal{R}\colon A \to A, \qquad
\mathcal{R}(x)=f_0^{-n(x)} f_1^{-m(x)}(x).
\end{equation*}
Note that this map can be written as
$
\mathcal{R}(x)=F^{n(x)+m(x)}(x)$ for $x\in A$,
where $F:[a,\,b]\to [a,\,b]$ is defined by $F=f_0^{-1}$ in $[a,\,f_1(a)]$ and $F=f_1^{-1}$ in $(f_1(a),\,b]$. Therefore, for every $x\in [a,\,b]$ there
is a smallest non-negative number $k(x) \geq 0$ such that $F^{k(x)}(x)\in A$ and $\mathcal{R}$ can be extended to the whole interval by taking
$$
   \mathcal{R}\colon(a,b)\to A, \qquad  \mathcal{R}(x)=F^{n+m+k}(x),
$$
where $k=k(x)$, $m=m(F^k(x))$ and $n=n(F^{k}(x))$.

A point $d\in A$ is said to be a \emph{discontinuity of
$\mathcal{R}$} if $\mathcal{R}(d)=f_0^{-1}f_1(a)$ or equivalently, if
$d=f_1^{m(d)}f_0^{n(d)-1}f_1(a)$. These points define a partition on $A$.
In other to describe this partition we have to consider two cases: $f_1^2(a)\not\in A$ and $f_1^2(a)\in A$. In the first case $m(x)=1$ for all $x\in A$ and we write $I_1=A$.
In the second case,
consider $\ell \in \mathbb{N}$ such that $f_1^{\ell+1}(a)\in A$, but $f_1^{\ell+2}(a)\not \in A$.
Then $f_1^mf_1(a)$ for $m=1,2,\ldots,\ell$ defines a partition on $A$ given by
\begin{align*}
I_\ell&=(f_1^{\ell+1}(a),\,f_0^{-1}f_1(a)] \quad
 \text{and} \\ I_{m}&=(f_1^{m}(a),\,f_1^{m+1}(a)] \quad \text{for} \ 1\leq m< \ell.
\end{align*}
Hence, $m(x)=m$ for each $x\in I_{m}$ and $f_1^{-m(x)}(I_{m}) \subset (a,\,f_1(a)]$.
Observe that the sequence of points $f_0^nf_1(a)$, $n\geq 0$, creates a partition on the interval $(a,\,f_1(a)]$.
This partition induces a partition in right-closed intervals $I_{mn}$ in each $I_m$. Namely,
$I_{mn} \subset f_1^{m}f_0^n(A)$  for $1\leq m \leq \ell$ and
$n\geq 1$.

\begin{prop}
\label{l-navas}
For every $J=[x,y] \subset I_{mn}$ it holds that
\begin{equation*}
  \frac{|\mathcal{R}(J)|}{|J|} \geq (1-\varepsilon) \varepsilon^{-1}  e^{-V} \, \frac{D(\mathcal{R}(y))}{D(y)}
\end{equation*}
where $D: A \to (0,\infty)$ is given by $D(y)=|f_0(y)-y|.$
\end{prop}
\begin{proof} Notice that $\mathcal{R}(J)=f_0^{-n}f_1^{-m}(J)$ and
hence
\begin{equation}
\label{eq:factors}
   \frac{|\mathcal{R}(J)|}{|J|}= \frac{|f_0^{-n}f_1^{-m}(J)|}{|f_1^{-m}(J)|} \cdot \frac{|f_1^{-m}(J)|}{|J|}.
\end{equation}
In order to estimate the first factor, observe that $f_1^{-m}(J)\subset [f_0f_1^{-m}(y),\, f_1^{-m}(y)]$ and so
\begin{align*}
 \mathrm{Dist}(f_0^{-n},&[f_0f_1^{-m}(y),\, f_1^{-m}(y)]) =  \mathrm{Dist}(f_0^n,f_0^{-n}([f_0f_1^{-m}(y),\, f_1^{-m}(y)]))
 \\ &\leq \sum_{i=0}^{n-1} \mathrm{Dist}(f_0,f_0^{i-n}([f_0f_1^{-m}(y),\, f_1^{-m}(y)]))
\leq V_0 |a-f_0^{-1}f_1(a)|.
\end{align*}
Thus,
\begin{equation}
\label{eq:firstfactorA}
 \frac{|f_0^{-n}f_1^{-m}(J)|}{|f_0^{-n}f_0f_1^{-m}(y)-f_0^{-n}f_1^{-m}(y)|} \geq e^{-V_0|a-f_0^{-1}f_1(a)|} \frac{|f_1^{-m}(J)|}{|f_0f_1^{-m}(y)-f_1^{-m}(y)|}.
\end{equation}
Let us  denote by $t=f_1^{-m}(y)$. Using the mean value theorem, there is  some $\xi \in (a,t)$ such that
$
f_0(t)-a=f_0(t)-f_0(a)=Df_0(\xi)(t-a)
$.
Hence
$$
f_0(t)-t = \frac{Df_0(\xi)-1}{Df_0(\xi)} \, (f_0(t)-a) \leq \frac{\varepsilon}{1-\varepsilon} \, (f_0(t)-a).
$$
Substituting in~\eqref{eq:firstfactorA} one has that
\begin{equation}
\label{eq:firstfactor}
 \frac{|f_0^{-n}f_1^{-m}(J)|}{|f_1^{-m}(J)|} \geq (1-\varepsilon)\varepsilon^{-1} e^{-V_0|a-f_0^{-1}f_1(a)|}\, \frac{|f_0(\mathcal{R}(y))-\mathcal{R}(y)|}{|f_0f_1^{-m}(y)-a|}.
\end{equation}
Now we will estimate the second factor in~\eqref{eq:factors}. Observe that $f_1^{-m}(J)$ and $[a,\,f_1^{-m}(y)]$ are contained in the fundamental domain  $[a,\,f_1(a)]$ of $f_1$. Thus,
\begin{align*}
 \mathrm{Dist}(f_1^{m},[a,\,f_1(a)]) &\leq   \sum_{i=0}^{m-1} \mathrm{Dist}(f_1, f_1^i([a,\,f_1(a)]))
\leq V_1  |a-f_0^{-1}f_1(a)|
\end{align*}
and so
\begin{equation*}
 \frac{|J|}{|f_1^m(a)-y|} = \frac{|f_1^{m}f_1^{-m}(J)|}{|f_1^m([a,\,f_1^{-m}(y)])|} \leq e^{V_1 |a-f_0^{-1}f_1(a)|} \, \frac{|f_1^{-m}(J)|}{|f_1^{-m}(y)-a|}.
\end{equation*}
As $a<f_0f_1^{-m}(y)<f_1^{-m}(y)$ and since $y\in I_{mn}$, by the construction of the partition $I_{mn}$ of $I_{m}$, one has that
$f_0(y)\leq f_1(a)\leq f_1^m(a) < y$.  Then the above inequality implies
\begin{equation}
\label{eq:secondfactor}
   \frac{|f_1^{-m}(J)|}{|J|} \geq e^{-V_1 |a-f_0^{-1}f_1(a)|} \, \frac{|f_0f_1^{-m}(y)-a|}{|f_0(y)-y|}
\end{equation}
Putting together~\eqref{eq:firstfactor} and
\eqref{eq:secondfactor} in Equation~\eqref{eq:factors} one deduces
the inequality desired and concludes the proof of the lemma.
\end{proof}

The following lemma shows that some iterate of the first return map $\mathcal{R}$ is expanding.

\begin{lemma}
\label{r:expansividad_R}
There exists $N \in \mathbb{N}$ such that  $D\mathcal{R}^N(x) >1$ for all $x \in A$.
\end{lemma}
\begin{proof} Let us denote by $C(f_0)= \inf \{ D(x): x\in A \} /
\sup \{ D(x) :  x \in A\}  >0$.  Fix $\kappa>1$ and let $N=N(f_0)$
be a natural number such that
$$
((1-\varepsilon)\varepsilon^{-1} e^{-V})^{N} \, C(f_0) >\kappa.
$$
Suppose $x\in I_{mn}$ and consider $\delta>0$ small enough such that $J=[x-\delta,x]\subset I_{mn}$. Hence, by Proposition~\ref{l-navas} it follows that
\begin{align*}
\frac{|\mathcal{R}^N(J)|}{|J|} &=
\frac{|\mathcal{R}^{N}(J)|}{|\mathcal{R}^{N-1}(J)|}\cdots  \frac{|\mathcal{R}(J)|}{|J|} \\ &\geq ((1-\varepsilon)\varepsilon^{-1} e^{-V})^{N} \, \frac{D(\mathcal{R}^N(x))}{D(x)}  \\
&\geq ((1-\varepsilon)\varepsilon^{-1} e^{-V})^{N} \, C(f_0) >\kappa.
\end{align*}
Since this inequality holds for every $\delta> 0$ small enough,
this allows to conclude that $D\mathcal{R}^N(x)\geq \kappa >1$.
\end{proof}

The first simplification to prove the minimality and the density of periodic points in a $**$-interval is to note that it
suffices to show that the orbit of the global attractor $a$ for $f_0|_{K^{**}}$,
with respect to $\IFS(\Phi)$ is dense in the interval $K^{**}=[a,b]$.

\begin{lemma}
\label{lem-densidad}
If $K^{**} \subset \overline{\mathrm{Orb}^+_{\Phi}(a)}$ then
\begin{equation*} 
  K^{**} \subset\overline{\mathrm{Per}(\IFS(\Phi))}
  \quad \text{and} \quad
  K^{**} \subset \overline{\mathrm{Orb}^+_{\Phi}(x)} \quad \text{for all } x \in K^{**}.
\end{equation*}
\end{lemma}
\begin{proof} Consider $x\in K^{**}$ and let $U$ be any open set
in $K^{**}$. From the density of the orbit of $a$ for $\IFS(\Phi)$
in $K^{**}$, there is $h\in \IFS(\Phi)$ such that $h(a) \in U$.
Since $h$ is a continuous map and $a$ is a global attractor of
$f_0$ in $K^{**}$, there exists $\ell\in \mathbb{N}$ such that
$f_0^\ell(x)$ is close enough to $0$, so that $h\circ f_0^\ell(x)
\in U$. Therefore, the orbit of $x$ for $\IFS(\Phi)$ is dense in
$K^{**}$.

Now, given $x\in K^{**}$ we will show that $ x \in
\overline{\mathrm{Per}(\IFS(\Phi))}$. Let $I$ be any open interval
such that $x\in I$. From the density of the orbit of $a$, there is
$h\in \IFS(\Phi)$ such that $h(a)\in I$. Since $h$ is a continuous
map there is $\delta>0$ such that $h((a-\delta,\,a+\delta))
\subset I$.  Using that $a$ is a global attractor of $f_0$, there
is $\ell\geq 0$ such that $f_0^\ell(I) \subset
(a-\delta,\,a+\delta)$. Then the closure of $h\circ f_0^\ell(I)$
is strictly contained in $I$. By the intermediate value theorem,
$h\circ f_0^\ell$ has a fixed point in $I$ and thus $I\cap
\mathrm{Per}(\IFS(\Phi))\not = \emptyset$. This implies that $x$
belongs to the closure of the periodic points of $\IFS(\Phi)$ and
we conclude the proof of the lemma. \end{proof}

Now, we are ready to proof Theorem~\ref{t.duminy_lemma}:

\begin{proof}[Proof of Theorem~\ref{t.duminy_lemma}] Recall that the
return map $\mathcal{R}:A\to A$ can be extended to the interval
$(a,b)$. In particular, this implies that for any open interval $I
\subset K^{**}$, there exists $h\in \IFS(\Phi)$  such that
$h^{-1}(I)\cap A\not = \emptyset$. From
Lemma~\ref{r:expansividad_R}, the return map $\mathcal{R}^N$ is an
expanding map in $A$. Thus, there is $\ell\in\mathbb{N}$ such that
$\mathcal{R}^{\ell N}(h^{-1}(I)\cap A)$  contains some
discontinuity of $\mathcal{R}^N$.

Recall that the discontinuities $d\in A$ are points in the orbit
of $a$ for $\IFS(\Phi)$, i.e., applying the appropriate inverse
branch,  $d=\mathcal{R}^{-N}(f_0^{-1}f_1(a))=L\circ f_1(a)$ for
some $L\in \IFS(\Phi)$. Then, one has $h\circ L \circ f_1(a) \in
I$. Therefore, the orbit of $a$ is dense in $K^{**}$ with respect
to $\IFS(\Phi)$. Finally, from Lemma~\ref{lem-densidad} we
conclude the proof of the theorem. \end{proof}

\section{Cycles for IFS and proof of Theorem~\ref{thmA}}
\label{ss:cycle}

Now we will give a proof of Theorem~\ref{thmA}. Observe that as an
$ss$-interval is a forward invariant set with respect to
$\IFS(f_0^{n_0},f_1^{n_1})$. Thus, the existence of an
$ss$-interval immediately implies that
$\IFS(f_0^{n_0},f_1^{n_1})$ is not minimal. From now on assume
that there are no $ss$-intervals for $\IFS(f_0^{n_0},f_1^{n_1})$,
and so we have to show minimality of~$\mathbb{S}^1$.

We will first show the result when the rotation numbers of $f_0$ and $f_1$ are zero,
that is, both of the diffeomorphisms have fixed points.
The proof will require the following notion of a cycle for an IFS
of two diffeomorphisms on the circle. Similarly to Duminy's result (Theorem~\ref{t.duminy_lemma}),
using the cycle we will then construct a return map and afterwards estimate the derivative.
Here the expanding return map is of a global character, while Duminy's lemma can be thought
of as a local version.

\subsection{Cycles.}
Since $f_0$ and  $f_1$ do not have fixed points in common, we can (without ambiguity)
denote by $\{s_i\}$ the set of fixed points of $f_0$ and $f_1$ with a non-empty basin of attraction.
This set includes the fixed points which may attract only from one-side (semi-attracting). Let $B(s_i)$
represent the basin of attraction of $s_i$ with respect to the relevant map $f_0$ or $f_1$ which does not include the point $s_i$ (it is an open set).

\begin{defi}[Cycle for IFS]
Define a partial order on the
(semi)attracting fixed points by $s_i \prec s_j$ if and only if $s_i$ belongs
to $B(s_j)$. A sequence of (semi)attractors is said to be a \emph{cycle
of length $n$ for $\mathrm{IFS}(\Phi)$} if
 $$s_{i_1}
\prec s_{i_2} \prec \cdots \prec s_{i_{n}} \prec
s_{i_{n+1}} \ \ \text{and} \ \ s_{i_1}=s_{i_{n+1}}.$$
\end{defi}

Notice that:
\begin{itemize}
 \item since the (semi)attracting points $s_i$ of the cycle alternate
 with respect to the maps $f_0$ and $f_1$, the length of a cycle is
 always even.
 \item a $ss$-interval is a cycle of length 2.
\end{itemize}

The following result shows the existence of a cycle under the assumption of no periodic points in common.

\begin{prop} \label{p.ciclo0}
Let $f_0$, $f_1$ be circle diffeomorphisms with zero rotation number and no fixed points
in common. Then there exists at least one cycle for
$\IFS(\Phi)$.
\end{prop}
\begin{proof} Since $f_0$ and $f_1$ have no fixed points in
common, given $x\in \mathbb{S}^1$, then $x\in B(s_{i})$ for some
$s_{i}$. By compactness we can take a finite sub-covering and
re-ordering the indices, we may suppose that $\mathbb{S}^1=
B(s_1)\cup\dots\cup B(s_n)$. For any $s_i$, there exists
$s_{i_2}\in \{s_1, \dots, s_n\}$ with $s_i\prec s_{i_2}$.
Proceeding inductively we obtain the sequence ($s_i=s_{i_1}$),
$s_{i_1} \prec s_{i_2} \prec \cdots \prec s_{i_{n}} \prec
s_{i_{n+1}}$. Then necessarily $s_{i_{n+1}}=s_{i_j}$ with $j\in
\{1,\dots,n\}$, and so we have
the cycle, 
$s_{i_j} \prec s_{i_2} \prec \cdots \prec s_{i_{n}} \prec
s_{i_j}$.
\end{proof}

By the above proposition we may assume that there exists a cycle
for $\mathrm{IFS}(\Phi)$. Let $n$ be the length of the cycle;
after re-numbering the indices $i_j$ of the (semi)attractors if
necessary, the cycle can be written as $ s_{n} \prec s_{n-1} \prec
\dots \prec s_2 \prec s_1 \prec s_0$ with $s_{0}=s_{n}$, where
$s_0$ is a (semi)attracting fixed point of $f_0$. Note that $s_k$
is a (semi)attracting fixed point for $f_{k\, \mathrm{mod}\, 2}$:
if $k$ is an even number then $s_k$ is a fixed point of $f_0$ and
if $k$ is an odd number then of
$f_1$.\\[-0.3cm]

\noindent \textbf{Notation:} For the rest of the section $f_k$
will stand for $f_{k\, \mathrm{mod}\, 2}$.  \medbreak

We consider $\mathbb{S}^1$ parameterized by
$[s_0,\,s_0+1]\,\mathrm{mod}\, 1$, and so $s_0 < s_k <s_{0}+1=s_n$
for $k=1,\ldots, n-1$ with respect to the real order on the
interval $[s_0,\,s_0+1]$. We denote by $s_k^{-}$ and $s_k^{+}$
(when they exist) the fixed points of the lift of $f_{k}$ on
$[s_0, s_0+1]$ closest to $s_k$ from the left and right,
respectively. If $f_k$ has more than one fixed point and since
$s_k$ is a (semi)attractor, necessarily one of the points,
$s_k^{-}$ or $s_k^{+}$, exists and is a (semi)repeller. With
respect to the order on the interval $[s_0,\, s_0+1]$, two (non
mutually excluding) cases are possible: $ s_1\in (s_0,\, s_0^{+})$
or $ s_1\in (s_n^-,\, s_n)$. We will assume the former case, the
latter being analogous. The next lemma on the order of the points
is crucial for the creation of the return map.

\begin{lemma}\label{claim:sk}
If there are no $ss$-intervals for $\mathrm{IFS}(\Phi)$ then
on the interval $[s_0,\, s_0+1]$ the cycle has the order
$$
s_0< s_1<\dots< s_{n-1}<s_0+1=s_n, \quad \text{where $s_{k+1}\prec s_k$ \ for \ $k=0,\dots,n-1$.}
$$
Moreover,
$$
 s_k <f_{k+1}^{-1}(s_k) < s_{k+1} < s_k^+ \quad \text{for} \ k=0,\ldots,n-1.
$$
\end{lemma}
\begin{proof} We will show that $s_{k+1} \in
(s_k,\,s_{k}^+)\subset B(s_k)$ for every $k=0,\dots,n-1$ which
proves the first part of the lemma. The second part will follow as
the fixed point of $f_{k+1}$ closest to $s_k$ from the right has
to be a repeller from the left (otherwise a $ss$-interval is
created) and so $s_k<f_{k+1}^{-1}(s_k)< s_{k+1}<s^+_k$.

First, we will assume that $f_0$ or $f_1$ has a unique fixed point. In this case, the cycle has length~$2$ and since there are no $ss$-intervals then $s_2=s_0+1$.
Thus $s_0<s_1<s_0^+\leq s_0+1=s_2$, where $s_2\prec s_1 \prec s_0$.
In particular, $s_{1} \in (s_0,\,s_{0}^+)\subset B(s_0)$ and we have $s_{2} \in (s_1,\,s_{1}^+)\subset B(s_1)$.

Now, we will assume that $f_0$ and $f_1$ have more than one fixed point. By the initial hypothesis we have
$s_1\in (s_0,\, s_0^+)\subset B(s_0)$. Let us show that
$s_2\in (s_1,\, s_1^+)\subset B(s_1)$.
Since $s_2\prec s_1$, $s_2\in B(s_1)\subset (s_1^-,\, s_1)\cup (s_1,\, s_1^+)$ and so
either $s_2\in (s_1^-,\, s_1)$ or $s_2\in (s_1,\, s_1^+)$.
It is enough to prove that the first case cannot occur.

Suppose that $s_2\in (s_1^-,\, s_1)$.
If $s_0 <s_1^-$ then $s_0 < s_2 <s_1 < s_0^+$, which is a contradiction since $s_0^+$ is the closest
fixed point of $f_0$ to $s_0$ from the right. If $s_0>s_1^-$, a similar argument concludes that $s_0=s_2$.
Since there are no $ss$-intervals and $s_2\prec s_1$ then $B(s_1)=(s_1,\,s^+_1)$ and $s_0+1 = s_2 +1 <s_1^+$.
Thus $s_1^-<s_0 < s_0+1 <s_1^+$, which implies that $f_1$ has a unique fixed point, a contradiction.
Therefore $s_2$ does not belong to $(s_1^-,\,s_1)$ and so it is contained in $(s_1,\,s^+_1)$.

Suppose by the inductive hypothesis that $s_i\in
(s_{i-1},\,s_{i-1}^+)\subset B(s_{i-1})$ for all $i\leq k$. To
prove that $s_{k+1}\in (s_{k},\,s_{k}^+)\subset B(s_{k})$ the idea
is basically the same as above. If $s_{k+1}\in (s_k^-,\, s_k)$,
then either $s_k^-\in(s_{k-1},\,s_k)$ or $s_k^-< s_{k-1}$. The
first case cannot occur because $s_{k+1}$ and $s_{k-1}$ are fixed
points of the same map $f_{k+1}=f_{k-1}$ and $s_{k+1}\in (s_k^-,\,
s_k)$. Thus $s_k^-< s_{k-1}$ and so $s_{k+1}\leq s_{k-1}$. This
fact together with $s_{k+1} \prec s_k$ imply that $(s_k^-,\,s_k)
\subset B(s_k)$ and hence $[s_{k-1}, \, s_k]$ is a $ss$-interval,
and thus we get a contradiction. Therefore $s_k<s_{k+1}<s_k^+$ and
so $(s_k,\,s_k^+)\subset B(s_k)$, which completes the induction
and concludes the proof of the lemma. \end{proof}

\subsection{Creating a return map.}
In the previous section was shown the existence of a cycle that
can be ordered in the following manner: $ s_{n} \prec s_{n-1}
\prec \dots \prec s_2 \prec s_1 \prec s_0$. Let us now use this
cycle to create a return map. Remember, $f_k$ stands for $f_{k\,
\mathrm{mod}\, 2}$.

\begin{lemma}\label{lema.creating_returning_map}
There exist families of right-closed pairwise disjoint intervals
$I_{i_1\ldots i_n} \subset \mathbb{S}^1$ and maps $h_{i_i \ldots
i_n} \in \mathrm{IFS}(\Phi)$ with $i_j \geq 0$ for $j=1,\ldots,n$
such that
\begin{enumerate}
\item  $A=(s_0,\,f_1^{-1}(s_0)]= \bigcup I_{i_1 \dots i_n}$,
\item $h_{i_1 \dots i_{n}}^{-1}(I_{i_1 \dots i_{n}}) \subset A$
if $i_{n}=0$ and $h_{i_1 \dots i_{n}}^{-1}(I_{i_1 \dots i_{n}})=A$
if $i_{n} > 0$,
\item if $c\in A\setminus \{ f_1^{-1}(s_0) \}$
is an endpoint of $I_{i_1 \dots i_{n}}$ then there exist a map $h \in
\mathrm{IFS}(\Phi)$ and a point $s \in \mathrm{Orb}_\Phi^+(s_1)\cup\dots\cup\mathrm{Orb}_\Phi^+(s_n)$ such that
$h(s)=c$. That is, the point $c$ is in the orbit of the
cycle for $\IFS(\Phi)$.
\end{enumerate}
\end{lemma}

\begin{proof}
We can define a non-empty fundamental domain of
$f_{k+1}$,
$$
A_k=(s_k,\,f^{-1}_{k+1}(s_k)] \quad \text{for $k=0,\dots,n$.}
$$
Since on the circle $s_0=s_n$, we identify on $\mathbb{S}^1$ the
intervals in the real line $A_0 =(s_0, f_1^{-1}(s_0)]$ and
$A_n=(s_n,f_1^{-1}(s_n)]$, denoting this interval by $A$.  In
order to create the expanding return map we will divide this
fundamental domain inductively.

By Lemma~\ref{claim:sk}, $s_0<f_1^{-1}(s_0)<s_1$ and since $s_1 \prec s_0$, there exists $j\in\mathbb{N}$ such that $
  s_0 < f_0^{j}(s_1)<f_1^{-1}(s_0)\leq f_0^{j-1}(s_1)<s_1.
$ Then
$$
    A=(s_0, \, f_{1}^{-1}(s_0)]=
    \bigcup_{i_1=0}^\infty I_{i_1},
$$
with
\begin{align*}
I_0 = (f_0^{j}(s_1),\, f_1^{-1}(s_0)]  \quad \text{and} \quad
I_{i_1} =
(f_0^{j+i_1}(s_1),\, f_0^{j+i_1-1}(s_1)] \ \ \text{if $i_1>0$}.
\end{align*}
Let $h_{i_1}=
f^{j+i_1}_0$. We have that $h^{-1}_{i_1}(I_{i_1})=(s_1,\,c_{i_1}]$
where
\begin{align*}
c_0 = f_0^{-{j}} f_1^{-1}(s_0) \in (s_1,\, f_0^{-1}(s_1)] \quad \text{and} \quad
c_{i_1}= f_0^{-1}(s_1) \ \ \text{if $i_1>0$.}
\end{align*}
Therefore, $h^{-1}_0(I_0)\subset A_1$ and $h^{-1}_{i_1}(I_1)=A_1$
if $i_1>0$.

Let $c \in A\setminus\{f^{-1}_1(s_0)\}$ be an endpoint
of $I_{i_1}$. Then, either, $c=f_0^{j}(s_1)$ if $i_1=0$ or
$c\in\{f_0^{j+i_1}(s_1),f_0^{j+i_1-1}(s_1)\}$ if $i_1>0$. In any
case, $c$ belongs to the orbit of the cycle.
This completes the first step of the induction and now we proceed
with the inductive hypothesis.

Suppose that we have families of right-closed pairwise disjoint
intervals $I_{i_1\ldots i_k} \subset \mathbb{S}^1$ and maps
$h_{i_1\ldots i_k}\in \IFS(\Phi)$ with $i_j \geq 0$ for
$j=1,\ldots,k$ such that
\begin{itemize}
\item[(i)] $A = \bigcup I_{i_1\dots i_k}$.
\item[(ii)]$h_{i_1 \dots i_{k}}^{-1}
(I_{i_1 \ldots i_{k}})= (s_k,\, c_{i_1\ldots i_{k}}]$ where
$$
c_{i_1\ldots i_{k-1}0} \in (s_k,\,f_{k+1}^{-1}(s_k)] \quad \text{and}
\quad c_{i_1\ldots i_{k}}=f_{k+1}^{-1}(s_k) \ \ \text{if $i_k>0$.}
$$In particular, $h_{i_1 \dots
i_{k-1}0}^{-1} (I_{i_1 \ldots i_{k-1}0})\subset A_k$ and $h_{i_1
\dots i_{k}}^{-1} (I_{i_1 \ldots i_{k}})=A_k$ if $i_k>0$.
\item[(iii)] If $c\in A\setminus \{ f_1^{-1}(s_{0})\}$
is an endpoint of $I_{i_1 \ldots i_{k}}$ then there exist a map $h \in
\IFS(\Phi)$ and a point  $s \in \mathrm{Orb}_\Phi^+(s_0)\cup\dots\cup \mathrm{Orb}_\Phi^+(s_n)$ such that
$h(s)=c$.
\end{itemize}
\noindent By Lemma~\ref{claim:sk} $s_k <
f^{-1}_{k+1}(s_{k})<s_{k+1}<s_k^+$. Hence, from the
inductive hypothesis we get
$$ s_k < c_{i_1 \ldots i_k} \leq
f^{-1}_{k+1}(s_{k})<s_{k+1}<s_k^+.$$
Now, since
$s_{k+1}\prec s_k$, for each multi-index $i_1 \ldots i_k$ there
exists  $j_{i_1\dots i_k} \in \mathbb{N}$ such that
$$
s_k < f^{j_{i_1\ldots i_k}}_{k}(s_{k+1})<
c_{i_1\ldots i_k} \leq f^{j_{i_1\ldots
i_{k}}-1}_{k}(s_{k+1})<s_{k+1}.
$$
Then
$$
h^{-1}_{i_1\ldots i_k}(I_{i_1\ldots i_k})=(s_k,\, c_{i_1 \ldots
i_k}]=\bigcup_{\ell=0}^\infty J_{i_1\dots i_k\ell}
$$
with
\begin{align*}
J_{i_1\dots i_k 0} &= (f^{j_{i_1 \dots
i_k}}_{k}(s_{k+1}), \,c_{i_1 \dots i_k}] \ \ \text{and} \\
J_{i_1 \ldots i_k \ell} &= (f_{k}^{j_{i_1 \dots
i_k}+\ell}(s_{k+1}),\, f_{k}^{j_{i_1 \ldots
i_k}+(\ell-1)}(s_{k+1})]
\ \ \text{if} \ \ell>0.
\end{align*}
Notice that the intervals $J_{i_1\ldots i_k\ell}$, $\ell
\geq 0$ are pairwise disjoint. Define $I_{i_1 \ldots i_k \ell}
\subset I_{i_1 \ldots i_k}$ by
$$
I_{i_1\ldots i_k \ell}=h_{i_1
\ldots i_k}(J_{i_1 \ldots i_k \ell}).
$$
By definition, for a fixed
multi-index $i_1\ldots i_k$, the intervals $I_{i_i\ldots i_k\ell}$,
$\ell\geq 0$ are also pairwise disjoint. Hence, by means of the induction
hypothesis on the intervals $I_{i_1\ldots i_k}$
and since $I_{i_1 \dots i_k \ell} \subset I_{i_1 \ldots i_k}$ it follows that
$
\{I_{i_1 \ldots i_{k+1}}: \ i_j \geq 0 \ \text{for} \
j=1,\ldots, k+1 \}
$
is a family of right-closed pairwise disjoint
intervals. Note that each right-closed interval $I_{i_1\ldots i_k}$
is the union of the intervals $I_{i_1\ldots i_k \ell}$, $\ell \geq
0$.
Then, by the induction hypothesis it also follows that $ A=\cup
I_{ i_1 \ldots i_{k+1}}$.

Set $h_{i_1\dots i_k i_{k+1}}=h_{i_1 \dots i_k} \circ
f_{k}^{j_{i_1 \dots i_k}+i_{k+1}}$. By
construction
$$
h_{i_1\dots i_{k+1}}^{-1} (I_{i_1, \dots i_{k+1}})=
f_{k}^{-j_{i_1 \dots i_k}-i_{k+1}} (J_{i_1\ldots
i_k i_{k+1}}) = (s_{k+1},\, c_{i_1 \dots i_{k+1}}]
$$
where $c_{i_1 \ldots i_k
i_{k+1}}=f_{k}^{-1}(s_{k+1})$ if $i_{k+1}>0$ and
$$
c_{i_1 \dots i_{k}0} = h^{-1}_{i_1\ldots i_{k}0} \circ
h_{i_1\ldots i_{k}}(c_{i_1\ldots i_k})=f^{-j_{i_1\ldots
i_k}}_{k}(c_{i_1\ldots i_k}) \in
(s_{k+1},\,f_{k}^{-1}(s_{k+1})].$$ Therefore,
$h_{i_1 \dots i_{k}0}^{-1} (I_{i_1 \ldots i_{k}0})\subset
A_{k+1}$ and $h_{i_1 \dots i_{k+1}}^{-1} (I_{i_1 \ldots
i_{k+1}})=A_{k+1}$ if $i_{k+1}>0$.

In order to prove the third item, we fix an interval $I_{i_1\ldots i_k\ell}$.
For every $\ell \geq 0$
the left endpoint of this interval is $h_{i_1 \dots i_k} \circ
f_{k}^{j_{i_1 \dots i_k}+\ell}(s_{k+1})$ and so
it belongs to the orbit of the cycle. The right endpoint is
$$
h_{i_1 \dots i_k}(c_{i_1\ldots i_k}) \ \ \text{if} \ \ \ell =0
\quad \text{and} \quad
h_{i_1 \dots i_k}\circ f_{k}^{j_{i_1 \dots
i_k}+\ell-1}(s_{k+1}) \ \ \text{if} \ \ \ell > 0.
$$
Note that $h_{i_1
\dots i_k}(c_{i_1 \dots i_k})$ is a endpoint of $I_{i_1 \dots
i_k}$. Therefore, by the inductive hypothesis,  either $h_{i_1
\dots i_k}(c_{i_1 \dots i_k})=f_1^{-1}(s_0)$ or $h_{i_1 \dots
i_k}(c_{i_1 \dots i_k})= h(s)$ for some $h \in \mathrm{IFS}(\Phi)$
and $s \in \cup_{i=0}^n \mathrm{Orb}_\Phi^+(s_i)$.

Going through the $n$ steps of the cycle we conclude the lemma.
\end{proof}

From the above lemma we define the return map over
$A=(s_0,\,f^{-1}_1(s_0)]$ as
$$
  \mathcal{R}\colon A \to A, \qquad \mathcal{R}|_{I_{i_1\ldots i_n}}=h^{-1}_{i_1\ldots i_n}.
$$
The endpoints of the intervals $I_{i_1\ldots i_n}$ are called
\emph{discontinuities} of $\mathcal{R}$ and if $s$ is a discontinuity
then $s \in \mathrm{Orb}_\Phi^+(s_0)\cup \dots\cup \mathrm{Orb}_\Phi^+(s_n)$.
We summarize some of the points from the inductive process in the construction of the return map
in Lemma~\ref{lema.creating_returning_map}.
This notation will be required for the next section.

\begin{addendum} \label{Add}
For each $k=1,\ldots, n$,  there exists a family
$$
\{ (I_{i_1\ldots i_k},h_{i_1\ldots i_k}, m_{i_1\ldots i_k}): i_j
\geq 0 \  j=1,\ldots k \}
$$
with $I_{i_1\ldots i_k}$ pairwise disjoint right-closed intervals
of $\mathbb{S}^1$, $h_{i_1\ldots i_k}\in \IFS(\Phi)$ and $m_{i_1
\ldots i_k}$ natural numbers such that
$A=(s_0,\,f_1^{-1}(s_0)]=\cup I_{i_1\ldots i_k}$, and for
$k=1,\ldots,n-1$
\begin{enumerate}
\item $I_{i_1\ldots i_k i_{k+1}} \subset I_{i_1\ldots i_k}$ and $I_{i_1}$, $h_{i_1\dots i_k}(I_{i_1\dots i_{k+1}})$ are
respectively contained in a fundamental domain of $f_0$ and $f_k$;
\item  $m_{i_1}=j+i_1$,
$m_{i_1\ldots i_k}=j_{i_1\ldots i_k} + i_{k+1}$, and $h^{-1}_{i_1}=f_0^{-m_{i_1}}$,
$$h^{-1}_{i_1\ldots i_{k}i_{k+1}}= f_{k}^{-m_{i_1\ldots i_{k+1}}}\circ   h^{-1}_{i_1\ldots i_{k}}=
f_{k}^{-m_{i_1\ldots i_{k+1}}}f_{k-1}^{-m_{i_1\ldots i_{k}}} \cdots
f_{1}^{-m_{i_1 i_2}} f_0^{-m_{i_1}};$$
\item  $h^{-1}_{i_1\ldots i_{k-1}0}(I_{i_1\ldots i_{k-1}0})=(s_k, \,c_{i_1\ldots i_k}]$ and  $h^{-1}_{i_1\ldots i_{k}}(I_{i_1\ldots i_{k}})=A_k$ if $i_k >0$ with
    $$
    s_k < c_{i_1\ldots i_k} \leq f_{k+1}^{-1}(s_k) \quad \text{and} \quad
    A_k=(s_k,\, f_{k+1}^{-1}(s_k)];
    $$
\item $m_{i_1\ldots i_k+1}= m_{i_1\ldots i_k}+1$ and   $m_{i_1\ldots i_{k-1}0}\geq 1$ satisfies that
$$
  f^{m_{i_1\ldots i_{k-1}0}}_{k}(s_{k+1}) <c_{i_1\ldots c_{i_k}} \leq  f^{m_{i_1\ldots i_{k-1}0}-1}_{k} (s_{k+1}).
$$
\end{enumerate}
\end{addendum}

\subsection{Estimation of the derivative for the return map.}
We will prove the following proposition which can be compared to
Proposition~\ref{l-navas}.

\begin{prop}
\label{l-navas2}
For every $J=[x,y] \subset I_{i_1\dots i_n}$ it holds that
\begin{equation*}
  \frac{|\mathcal{R}(J)|}{|J|} \geq \bigg(\frac{e^{-(V_0+V_1)}}{e^{V_1}-1} \bigg)^{n/2} \cdot \frac{D(\mathcal{R}(y))}{D(y)}
\end{equation*}
where $D: A \to (0,\infty)$ is given by $D(y)=|f_0(y)-y|$, $n$ is the length of
the cycle (always even), and $V_0, V_1$ are the distortion
constants of $f_0, f_1$ respectively. 
\end{prop}

Define by $\mathcal{R}_j=h_{i_1\dots i_j}^{-1}$, $1\leq j\leq n$, and then $\mathcal{R}_n=\mathcal{R}$. We adopt the  convention $\mathcal{R}_0=\mathrm{id}$.  To obtain Proposition~\ref{l-navas2} we need the following two-by-two term estimate.

\begin{lemma}
 \label{l-navas3}
 For any $0\leq j \leq n-2$ even,
 $$\frac{|\mathcal{R}_{j+2}(J)|}{|\mathcal{R}_j(J)|}
 \geq \frac{e^{-(V_0+V_1)}}{e^{V_1}-1}\cdot \frac{D(\mathcal{R}_{j+2}(y))}{D(\mathcal{R}_j(y))}.$$
\end{lemma}

This lemma concludes immediately Proposition~\ref{l-navas2}. Indeed, we have that
$$\frac{|\mathcal{R}(J)|}{|J|}=\frac{|\mathcal{R}_{n}(J)|}{|\mathcal{R}_{n-2}(J)|}
\cdot \frac{|\mathcal{R}_{n-2}(J)|}{|\mathcal{R}_{n-4}(J)|}\cdots
\frac{|\mathcal{R}_{2}(J)|}{|J|}$$
and applying Lemma~\ref{l-navas3} to each of the terms one gets the estimate required in Proposition~\ref{l-navas2}.

\begin{proof}[Proof of Lemma~\ref{l-navas3}] Since
$$\frac{|\mathcal{R}_{j+2}(J)|}{|\mathcal{R}_j(J)|}=
\frac{|\mathcal{R}_{j+2}(J)|}{|\mathcal{R}_{j+1}(J)|}\cdot
\frac{|\mathcal{R}_{j+1}(J)|}{|\mathcal{R}_j(J)|},$$
the lemma will be proved by estimating each term.

Observe that since the cycle begins with the map $f_0$, necessarily $j+2$ (mod $2)=0$ and $j+1$ (mod $2)=1$.
Since $s_j$ is a fixed point of $f_0$ and $\mathcal{R}_j(J)$ is contained in
a fundamental domain of $f_0$ (see Addendum~\ref{Add}) then
$$s_j< f_0(\mathcal{R}_j(y))< \mathcal{R}_j(x)< \mathcal{R}_j(y).$$
Applying
$f_0^{-m_{i_1\dots i_{j+1}}}$ to the fundamental domain
$[f_0(\mathcal{R}_j(y)), \,\mathcal{R}_j(y)]$, the usual distortion estimate gives
$$\frac{| f_0(\mathcal{R}_j(y))-\mathcal{R}_j(y)|}{| \mathcal{R}_j(x)-\mathcal{R}_j(y)|} \geq
e^{-V_0}\frac{| f_0^{-m_{i_1\dots i_{j+1}}}(f_0(\mathcal{R}_j(y)))-f_0^{-m_{i_1\dots i_{j+1}}}(\mathcal{R}_j(y))|}
{| f_0^{-m_{i_1\dots i_{j+1}}}(\mathcal{R}_j(y))-f_0^{-m_{i_1\dots i_{j+1}}}(\mathcal{R}_j(x))|}.$$
Observe that
$$f_0^{-m_{i_1\dots i_{j+1}}+1}(\mathcal{R}_j(y))< s_{j+1}< f_0^{-m_{i_1\dots i_{j+1}}}(\mathcal{R}_j(y)) \quad \text{and} \quad
\mathcal{R}_{j+1}=f_0^{-m_{i_1\dots i_{j+1}}}\circ\mathcal{R}_j.
$$
Using this information and rearranging
the terms one obtains,
$$\frac{| \mathcal{R}_{j+1}(J)|}{|\mathcal{R}_j(J)|}\geq e^{-V_0}\frac{| s_{j+1}-\mathcal{R}_{j+1}(y)|}
{D(\mathcal{R}_j(y))}.$$
To estimate the other term, since $s_{j+2}$ is again a fixed point of $f_0$ and $\mathcal{R}_{j+2}(J)$ is
contained in a fundamental domain of $f_0$, it follows
$$s_{j+2}< f_0(\mathcal{R}_{j+2}(y))< \mathcal{R}_{j+2}(x)< \mathcal{R}_{j+2}(y).$$
By construction of the return map (see Addendum~\ref{Add}), $[s_{j+2}, \, \mathcal{R}_{j+2}(y)]\subset
A_{j+2}$, which is a fundamental domain of $f_1$.

This time let us apply
$f_1^{m_{i_1\dots i_{j+2}}}$ to $[f_0(\mathcal{R}_{j+2}(y)),\, \mathcal{R}_{j+2}(y)]$, contained
in a fundamental domain of $f_1$. Then the distortion estimate becomes
$$\frac{| \mathcal{R}_{j+2}(x)-\mathcal{R}_{j+2}(y)|}
{| f_0(\mathcal{R}_{j+2}(y))-\mathcal{R}_{j+2}(y)|}
\geq e^{-V_1}\frac
{| f_1^{m_{i_1\dots i_{j+2}}}(\mathcal{R}_{j+2}(y))-f_1^{m_{i_1\dots i_{j+2}}}(\mathcal{R}_{j+2}(x))|}
{| f_1^{m_{i_1\dots i_{j+2}}}(f_0(\mathcal{R}_{j+2}(y)))-f_1^{m_{i_1\dots i_{j+2}}}(\mathcal{R}_{j+2}(y))|}.
$$
Since $A_{j+2}$ is a fundamental domain of $f_1$,
$$f_1(\mathcal{R}_{j+2}(y))<f_0(\mathcal{R}_{j+2}(y)< \mathcal{R}_{j+2}(y).$$
Therefore,
\begin{align*}
| f_1^{m_{i_1\dots i_{j+2}}}(f_0(\mathcal{R}_{j+2}(y)))&-f_1^{m_{i_1\dots i_{j+2}}}(\mathcal{R}_{j+2}(y))|< \\
&<| f_1^{m_{i_1\dots i_{j+2}}}(f_1(\mathcal{R}_{j+2}(y)))-f_1^{m_{i_1\dots i_{j+2}}}(\mathcal{R}_{j+2}(y))|
\end{align*}
As $f_1^{m_{i_1\dots i_{j+2}}}\circ \mathcal{R}_{j+2}=\mathcal{R}_{j+1}$, using the above inequality
and rearranging the terms,
$$\frac{| \mathcal{R}_{j+2}(J)|}{|\mathcal{R}_{j+1}(J)|}\geq e^{-V_1}\frac{D(\mathcal{R}_{j+2}(y))}
{| f_1(\mathcal{R}_{j+1}(y))-\mathcal{R}_{j+1}(y))|}.$$
Thus,
 $$\frac{|\mathcal{R}_{j+2}(J)|}{|\mathcal{R}_j(J)|}\geq e^{-(V_0+V_1)}\cdot
\frac{D(\mathcal{R}_{j+2}(y))}{D(\mathcal{R}_j(y))}\cdot
 \frac{| s_{j+1}-\mathcal{R}_{j+1}(y)|}
 {| f_1(\mathcal{R}_{j+1}(y))-\mathcal{R}_{j+1}(y))|}.$$
Finally let us estimate the last term. Applying the mean-value theorem to the function
$f_1-\mathrm{id}$ and using that $s_{j+1}$ is a fixed point of $f_1$,
$$ \frac{| s_{j+1}-\mathcal{R}_{j+1}(y)|}
 {| f_1(\mathcal{R}_{j+1}(y))-\mathcal{R}_{j+1}(y))|}=
 \frac{1}{| Df_1(z)-1|}>\frac{1}{e^{V_1}-1},$$
 which proves the lemma.
\end{proof}

Assume the maps $f_0, f_1$ are close enough to rotations
in the $C^{1+\mathrm{bv}}$-topology so that
$$
\frac{e^{-(V_0+V_1)}}{e^{V_1}-1}>1.
$$
Notice that, if the distortion constants $V_0\leq \varepsilon ,V_1 \leq \varepsilon$, the above condition is fulfilled for every positive $\varepsilon\leq 0{.}38$.  Then also
$\big((e^{V_1}-1)e^{(V_0+V_1)}\big)^{-n/2}>1$.
By the same proof of Lemma~\ref{r:expansividad_R}, we conclude the following
analogous  result.
\begin{lemma}\label{expand2}
For any $\lambda>1$,
there exists $N \in \mathbb{N}$ such that  $D\mathcal{R}^N(x) >\lambda$ for all $x \in A$.
\end{lemma}

\subsection{Proof of Theorem~\ref{thmA}}
Let us continue assuming that the maps $f_i$ both have fixed points. We will need the next
two lemmas.

\begin{lemma}\label{ssuu}
 There is a $ss$-interval for $\IFS(\Phi$) if and only if there is
 an $uu$-interval.
\end{lemma}
\begin{proof}
 By contradiction, suppose there exists a $uu$-interval but no $ss$-intervals
 (the other case is similar). Then, Lemma~\ref{claim:sk} implies that there is a cycle of length $n$ with the order
 $s_0<s_1<\dots<s_0+1=s_n$ on the interval $[s_0,\, s_0+1]$ identified with the circle.

 Consider the maps $f_i$ as maps of the real line (in particular of the interval $[s_0,\, s_0+1]$), and the
 $uu$-interval as the interval $[a,\, b]\subset [s_0,\, s_0+1]$. By the definition of the $uu$-interval, if $x\leq a$, then $f_i(x)\leq a$ and if $x\geq b$, then
 $f_i(x)\geq b$.
And so if $x>b$, then any finite composition of the maps $f_0$ and $f_1$, denoted by $h$, will satisfy that $h(x)>b$.

 Take $x>s_{n-1}$ and $x>b$. Using the cycle we can inductively attract $x$ to $s_{n-2}$, then to $s_{n-3}$, until $s_0$,
 such that in the end will obtain a map $h$, a finite composition of $f_0$ and $f_1$, with $h(x)<a$. This contradicts the previous observation.
\end{proof}

\begin{lemma}\label{ssuu2}
Suppose that there are no $ss$-intervals for $\IFS(\Phi)$ and for
some interval $I$,  $f_i(I)\cap I\neq\emptyset$, $i=0, 1$. Then
for every $x\in \mathbb{S}^1$, there exists $h_1,
h_2\in\IFS(\Phi)$ such that $h_1(x)\in I$ and $h_2^{-1}(x)\in I$.
\end{lemma}
\begin{proof} Since there are no $ss$-intervals,
Lemma~\ref{claim:sk} implies that there is a cycle with the order
$s_0<s_1<\dots<s_n+1$ in the interval  $[s_0, \, s_0+1]$. We
identify the interval $I$ on the circle $\mathbb{S}^1$ with an
interval in $[s_0,\,s_0+1]$ that has the endpoints $a<b$.

Assume first that on the interval  $[s_0,\,s_0+1]$, the point $x$ satisfies $b< x\leq s_0+1$.
Considering $f_i$ as maps of the interval and
arguing as in the proof of the above lemma, by using the attractors of the cycle we can inductively
start attracting $x$ to $s_0=s_n$. Then,
there exists $h$ a finite composition of $f_0$ and $f_1$ of the form $h=f_i\circ h_1$ ($i=0$ or $1$),
such that $h(x)\geq a$ and $h(x)< a$. Since $f_i(I)\cap I\neq\emptyset$, necessarily $h_1(x)\in I$,
as required.

For the case when $x\leq a$, suppose $I\subset [s_{i-1}, s_i]$. By
means of the cycle but now looking on the circle $\mathbb{S}^1$,
we can attract $x$ so $s_i$ so that on the interval
$[s_0,\,s_0+1]$, $x>s_i$. Then this is similar to the situation
above, applied not to $x$ but to $h(x)$, with $h(x)>s_i\geq b$.

Since by Lemma~\ref{ssuu} having no $ss$-intervals is equivalent
to having no $uu$-intervals, the same argument can be repeated
with $\IFS(\Phi^{-1})$ to prove the existence of $h_2^{-1}$ with
$h_2^{-1}(x)\in I$. \end{proof}

\subsubsection{Proof of minimality.}
Let us show that $\IFS(\Phi)$ is minimal. The minimality of
$\IFS(\Phi^{-1})$ follows analogously since by Lemma~\ref{ssuu},
$\IFS(\Phi^{-1})$ has no $ss$-intervals.

Consider a point $x\in \mathbb{S}^1$ and an open interval $J$.
Since there is no $ss$-intervals, by Lemma~\ref{expand2}, there
exists an expanding return map $\mathcal{R}^N$ defined in the
domain $A$ and associated with a cycle. Since $A=(s_0, \,
f^{-1}(s_0)]$, a fundamental domain of $f_1$ and is in the basin
of attraction of $s_0$, then $f_i(A)\cap A\neq\emptyset$. By
Lemma~\ref{ssuu2}, there exists $h_1\in \IFS(\Phi)$ with
$h_1^{-1}(J)\cap A\neq\emptyset$. We can assume $J$ is small
enough so that $h_1^{-1}(J)$ is contained in one of the domains of
the map $\mathcal{R}^N$. Since $D\mathcal{R}^N>\lambda>1$, then
$|\mathcal{R}^{kN}(J)|>\lambda^k|J|$. Iterating sufficiently many
times $\mathcal{R}^{kN}(J)$ eventually must contain a
discontinuity $s$ of $\mathcal{R}^N$.

By Lemma~\ref{lema.creating_returning_map}, $s\in  \mathrm{Orb}_\Phi^+(s_0)\cup \dots\cup \mathrm{Orb}_\Phi^+(s_n)$. Without loss of generality we
assume that $s$ is in the orbit of $s_0$.
Then there exists $h_2\in\IFS(\Phi)$ with $h_2(s_0)=s$ and
so
$$
h_2(s_0)\in\mathcal{R}^{kN}\circ h_1^{-1}(J).
$$
Writing $\mathcal{R}^{kN}$ as $h_3^{-1}$ with $h_3\in\IFS(\Phi)$,
we get  $h_1\circ h_3\circ h_2(s_0)\in J$.

Being the order $s_0<s_1<\dots<s_0+1=s_n$,
given a point $x$ we may attract $x$ arbitrary close to $s_0$
by inductively attracting to each $s_i$ of the cycle. Then there
exists $h_4$ with $h_4(x)$ close enough to $s_0$. It can be
concluded that $h_4\circ h_1\circ h_3\circ h_2(x)\in J$, and since
$x$ and $J$ were arbitrary, this proves minimality of
$\IFS(\Phi)$.  \medbreak

\begin{rem}
\label{exp} The action of the semigroup $\IFS(\Phi)$ is expanding:
for any $x$ in the circle, there exists $g\in\IFS(\Phi)$ such
 that $Dg^{-1}(x)>1$.
\end{rem}

To prove this, let us use the fact $\mathcal{R}$ is expanding
return map on some interval $A$. The branches of this expanding
return map are maps in $\IFS(\Phi^{-1})$.
 Since  $\IFS(\Phi^{-1})$ is minimal, there is
 $h\in\IFS(\Phi)$ such that $h^{-1}(x)\in A$. Then, for $n$ sufficiently large,
 $\mathcal{R}^n\circ h^{-1}(x)$ will have derivative larger than one.

\subsubsection{Density of hyperbolic periodic points.}
\label{ss:density hyp-points}
With respect to the return map $\mathcal{R}$, defined on the interval $A=(s_0,\,f_1^{-1}(s_0)]$,
consider one of the branches $I_{i_1 \dots i_{n}}$, such that $h_{i_1 \dots i_{n}}^{-1}(I_{i_1 \dots i_{n}})=A$.
Since  $h_{i_1 \dots i_{n}}^{-1}$ is expanding, there exists exactly one repelling fixed point, say $s$, of
 $h_{i_1 \dots i_{n}}^{-1}$. Then $s$ is an attracting fixed point of  $h_{i_1 \dots i_{n}}$, whose basin
 of attraction includes $A$.

 Now given an interval $J$, by minimality of $\IFS(\Phi)$ and assuming that $J$ is small enough,
 there exists $h_1\in\IFS(\Phi)$ such that $h_1(J)\subset A$.
 Again by minimality of $\IFS(\Phi)$,
 there exists $h_2\in \IFS(\Phi)$
 such that $h_2(s)\in J$.

Since $h_1(J)$ is contained in the basin of attraction of $s$, we may assume
that $h_1(J)$ is attracted sufficiently close to $s$, so that
$h_2\circ h_1(J)\subset J$ and $D(h_2\circ h_1)|_J<\beta< 1$.
This implies that there is a hyperbolic attractor in $J$.
The density of hyperbolic repellers for $\IFS(\Phi)$ follows
from the density of hyperbolic attractors for $\IFS(\Phi^{-1})$.

\subsubsection{The case of generators with periodic points.}
When the generators $f_i$
have non-zero rotation number, consider
the next lemma whose demonstration is left to the reader.

\begin{lemma}
\label{lem:noprove}
Let $f$ be an orientation preserving $C^{1+\mathrm{bv}}$-diffeomorphism of the circle with periodic points of period $n$, and consider $B$
the immediate basin of attraction of some periodic point.
Then
$$
e^{-V_f} \leq \frac{Df^n(x)}{Df^n(y)} \leq e^{V_f} \quad \text{for all $x,y \in B$}
$$
where $V_f$ is the distortion of $f$ on $\mathbb{S}^1$. That is,
$\mathrm{Dist}(f^n,B)\leq V_f$.

Since in the immediate basin there is always a point with derivative one,
$$
e^{-V_f} \leq Df^n(x) \leq e^{V_f} \quad \text{for all $x\in B$}.
$$
\end{lemma}

The lemma states that the
the distortion of $f_i^{n_i}$ is the same of $f_i$ in the relevant regions and the derivative
of $f_i^{n_i}$ is bounded by the distortion and thus will still be close to one if $f_i$ is close
to rotations. Namely, $|Df_i^n(x)-1| \leq e^{V_i}-1$ where $V_i$ is the distortion constant of $f_i$. Therefore, the same proof as in the case of fixed points applies.
Observe that also in Lemma~\ref{l-navas2} we have to ask again
that $e^{-(V_0+V_1)}/{e^{V_1}-1}>1$.
Since $V_0\leq \varepsilon$,  $V_1\leq \varepsilon$, the above
condition is fulfilled for all $\varepsilon \leq 0{.}38$.

This completes the proof of Theorem~\ref{thmA}.

\section{Spectral decomposition}
\label{s:spectral-decomposition}

\subsection{Spectral decomposition on the real line.}

The next theorem  gives a complete description of the global topological dynamics of an IFS generated by a pair of diffeomorphisms on the real line.

\begin{theo}[Spectral decomposition on the real line]
\label{thm:dec-real-line} There is $\varepsilon>0{.}38$ such that
if $f_0$ and $f_1$ are diffeomorphisms on the real line, one of
them having at least a fixed point, with no fixed points in common
and $\varepsilon$-close to the identity in the
$C^{1+\mathrm{bv}}$-topology, then
$$
    L(\IFS(\Phi))\setminus \{\pm \infty\}= \bigcup K_i \quad \text{and} \quad K_i \subset \overline{\mathrm{Per}(\IFS(\Phi))}
$$
where each $K_i$ is a compact,
maximal transitive set for $\IFS(\Phi)$. Moreover, each compact set $K_i$ is either
\begin{itemize}
\item a single fixed point of $f_0$ or $f_1$ or
\item a ${**}$-interval for $\IFS(\Phi)$ with $**\in \{ss,su,uu, s, u\}$.
\end{itemize}
In particular if $f_0$ and $f_1$ have only hyperbolic fixed points,
then the above union of the sets $K_i$ is countable and they are pairwise disjoint.
%
\end{theo}
\begin{proof} Consider $z\in L(\IFS(\Phi))\setminus\{\pm\infty\}$.
We can assume that $z\in \omega(\IFS(\Phi))$ since the situation
for the $\alpha$-limit of $\IFS(\Phi)$ is settled by a similar
argument. Then, by definition of $\omega$-limit of $\IFS(\Phi)$,
the point $z$ is approximated by points of the form $y_{k} \in
\omega_\Phi(x_k)$. For each $k$, there is a sequence $(h_{kn})_n
\subset \IFS(\Phi)$ such that each point $y_k$ is again
approximated by points of the form $h_{kn}\circ \dots \circ
h_{k1}(x_k)$.

\begin{claim}
\label{claim:spectral}
If $y\in \omega_\Phi(x)$, then either $y$ belongs to some $**$-interval for $\IFS(\Phi)$ with $**\in \{ss, su, uu, s, u\}$ or it is a fixed point of $f_0$ or $f_1$.
\end{claim}

This claim concludes that the limit set of $\IFS(\Phi)$ is contained in the union of $**$-intervals and the set of fixed points of $f_0$ and $f_1$.
Indeed, from Claim~\ref{claim:spectral}, either, $z$ is a fixed point of $f_0$ or $f_1$, or then, for $k_0$ large enough,
$y_k$ belongs in the same $**$-interval for all $k\geq k_0$ and thus $z$ is also in this $**$-interval.

\begin{proof}[Proof of the Claim~\ref{claim:spectral}] Since $y\in
\omega_\Phi(x)$, there is $(h_n)_n \subset \IFS(\Phi)$ such that
$y=\lim_{n\to \infty} h_n\circ\dots \circ h_1(x)$.  If $y$ is an
accumulation point of fixed points then, obviously, $y$ is a fixed
point of $f_0$ or $f_1$. Thus, we can assume that $y$ is not of
this type.
Also, without loss of generality, we may suppose that $y\geq 0$.

Since at least one of the generators has a fixed point, we can assume by the following argument
that $y\in [p,\,q]$ where $p,q \in \mathrm{Fix}(f_0)\cup \mathrm{Fix}(f_1)$ and $(p,\,q)\cap \mathrm{Fix}(f_i)=\emptyset$ for $i=0,1$.
Otherwise, there is a fixed point $P$ such that any other fixed point is less then $P$  and $y\in [P,\,\infty)$.
If $y\not=P$, then it is not hard to check via the geometry of the functions, that
it is not possible for $f_0, f_1 > \mathrm{id}$ in $(P\, \infty)$, since then $y=\infty$.
In the other case, $[P,\,\infty)$ is $s$ or $u$-interval or $y=P$.

Suppose $p$ and $q$ are both attractors or repellers but for
different maps restricted to $[p,\,q]$. Being $f_0$ and $f_1$
close to the identity, this closed interval is a ${ss}$ or
${uu}$-interval. So, we can consider the final case when $p$ and
$q$ are an attractor-repeller pair for the same map, say $f_0$
(the repeller-atractor case is analogous). Note that then $f_0 <
\mathrm{id}$ in $(p,\,q)$. We have two options: either $f_1
<\mathrm{id}$ or $f_1 > \mathrm{id}$ in $[p,\,q]$. In the first
case, both maps are below the identity and hence this geometry
implies that the unique possible $\omega$-limit points for
$\IFS(\Phi)$ in $[p,\,q]$ are $p$ and $q$.
%
%
For the second case, 
we have again two options:
$$
f_1([p,\,q]) \cap [p,\,q]\not = \emptyset \quad \text{or} \quad  f_1([p,\,q]) \cap [p,\,q]=\emptyset.
$$
In the first option $[p,\,q]$ is a ${su}$-interval for $\IFS(\Phi)$, and for the other one it follows as before that $y$ is either $p$ or $q$.

Therefore we have showed that in every possible case, the point
$y\in\omega_\Phi(x)$ belongs to a ${**}$-interval for $\IFS(\Phi)$
with $**\in \{ss,su,uu,s,u\}$ or it is a fixed point of $f_0$ or
$f_1$. \end{proof}

To conclude the proof of the theorem we need to show every set $K_i$ is contained in the limit set of $\IFS(\Phi)$.
Recall that $K_i$ denotes both, a \mbox{$**$-interval} or a fixed point of $f_0$ or $f_1$.
It is clear that every fixed point is contained in the limit set and so we only need to prove
that every $**$-interval for $\IFS(\Phi)$ is recurrent.

We will show that there exists $\varepsilon>0$ such that if
$f_0$ and $f_1$ are $\varepsilon$-close to the identity in the $C^{1+\mathrm{bv}}$-topology,
then any $**$-interval for $\IFS(\Phi)$ with $**\in \{ss,su,s\}$ satisfies the hypothesis of Theorem~\ref{t.duminy_lemma}.
This proximity implies that $|f_i(x)-x|<\varepsilon$,  $|Df_i(x)-1|<\varepsilon$ and $\mathrm{Dist}(f_i,\mathbb{R})<\varepsilon$ for $i=0,1$,
and we only need to show that $(1-\varepsilon)\varepsilon^{-1}e^{-V}>1$.

Notice that
$V\leq 2\varepsilon |a-f_i^{-1}f_j(a)|$ where $i\not=j$ and the $a$ is the endpoint of the $**$-interval.
Suppose $a$ is the left endpoint and $f_0<\mathrm{id}$ in $(a,\,f_0^{-1}f_1(a)]$.
From the proximity to the identity, $f_1(a)\leq a +\varepsilon$. Since
$$
f_0(x) \geq (1-\varepsilon)(x-a)+a \quad \text{for all  $x\in (a,\,f_0^{-1}f_1(a)]$},
$$
it follows that $f_0^{-1}f_1(a)-a \leq \varepsilon(1-\varepsilon)^{-1}$. Thus,
$$
  (1-\varepsilon)\varepsilon^{-1}e^{-V}\geq (1-\varepsilon)\varepsilon^{-1} e^{-2\varepsilon^2(1-\varepsilon)^{-1}}>1
$$
for all $\varepsilon\leq 0{.}38$.

Finally by Theorem~\ref{t.duminy_lemma}, ${ss}$, ${su}$ and ${s}$-intervals are minimal with a
dense set of periodic points. In particular, this kind of intervals are transitive sets for $\IFS(\Phi)$.
The same properties  are obtained for ${uu}$-intervals and $u$-intervals for $\IFS(\Phi^{-1})$.

\begin{lemma}
If $A$ is minimal for $\IFS(\Phi^{-1})$ then $A$ is  transitive for $\IFS(\Phi)$.
\end{lemma}

This lemma follows from the facts that minimality implies transitivity and if $A$ is transitive
for $\IFS(\Phi^{-1})$ then $A$ is  transitive for $\IFS(\Phi)$.

Applying the above lemma we obtain that ${uu}$ and ${u}$-intervals
are also transitive sets for the $\IFS(\Phi)$. This basically
concludes the proof of the theorem. The last part is immediately
obtained since if all the fixed points are hyperbolic, then they
are isolated and so there exists only a countable number of the
fixed points and the $**$-intervals. \end{proof}

\begin{rem}
Theorem~\ref{thm:dec-real-line} holds also for IFS generated by a pair of diffeomorphisms of the compact interval $I$ under the same assumptions.
\end{rem}
To see this, we understand the compact interval $I$ like the
compactified real line $[-\infty,\,\infty]$ and the $K^{s}$ or
$K^u$ intervals in $I$  are identified with the relevant intervals
on the real line.

\subsection{Spectral decomposition on the circle: Proof of Theorem~\ref{thmB}.}

The following result is a version of Theorem~\ref{t.duminy_lemma} for diffeomorphisms of the circle with periodic points.

\begin{prop}
\label{t.duminy_lemma-per} There is $\varepsilon>0{.}30$ such that
if $f_0$ and $f_1$ are circle diffeomorphisms with periodic points
of period $n_0$ and $n_1$ respectively, $\varepsilon$-close to the
rotation in the $C^{1+\mathrm{bv}}$-topology then for every
$**$-interval $K^{**}$, $**\in \{ss,su\}$, with respect to
$\IFS(\Phi^n)$ it holds
\begin{equation*}
\label{eq:den-min-Phi}
K^{**}  \subset
\overline{\mathrm{Per}(\mathrm{IFS}(\Phi^n))} \quad \text{and} \quad
K^{**} \subset \overline{\mathrm{Orb^+_{\Phi^n}}(x)} \ \ \text{for all} \
x\in K^{**}.
\end{equation*}
\end{prop}
\begin{proof} We will apply Theorem~\ref{t.duminy_lemma} for
$\IFS(\Phi^n)$ for the $**$-interval $K^{**}=[a,\,b]$. Without
loss of generality, we assume that $f^{n_0}_0 < \mathrm{id}$ in
$(a,\,b)$. According to Lemma~\ref{lem:noprove},
$$
   e^{-V_0}\leq Df_0^{n_0}(x) \leq  e^{V_0} \quad \text{for all $x\in \mathbb{S}^1$}
$$
where $V_0$ is the distortion constant of $f_0$. This implies that
$|Df^{n_0}_0(x)-1|<e^{V_0}-1$ for all $x\in \mathbb{S}^1$.
Similarly we denote by $V_1$ the maximal variation of $f_1$ and
then it also follows
$$
\mathrm{Dist}(f^{n_0}_0,I)\leq V_0 \quad \text{and} \quad \mathrm{Dist}(f^{n_1}_1,J)\leq V_1
$$
for all fundamental domains $I$ of $f_0^{n_0}$ and $J$ of $f_1^{n_1}$ in $(a,\,f_0^{-n_0}f_1^{n_1}(a)]$.

Since $(a,\,f_0^{-n_0}f_1^{n_1}(a)] \subset \mathbb{S}^1$, its
length is less that one. Let us remember the necessary condition
in Theorem~\ref{t.duminy_lemma},
$$
     (1-(e^{V_0}-1))(e^{V_0}-1)^{-1} e^{-(V_0 + V_1)|a-f_0^{-n_0}f_1^{n_1}(a)|} > 1.
$$
As $V_0,V_1 \leq \varepsilon$, it suffices to take
$\varepsilon\leq 0.{30}$ to guarantee that the above condition is
fulfilled, which concludes the proof of the proposition.
\end{proof}

\begin{proof}[Proof of Theorem~\ref{thmB}] This result is immediately
achieved from Theorem~\ref{thm:dec-real-line},
Proposition~\ref{t.duminy_lemma-per} and Theorem~\ref{thmA}.
Indeed, consider $\tilde{f}_0^{n_0}$ and $\tilde{f}_1^{n_1}$ as
the lifts on the real line of $f_0^{n_0}$ and $f_1^{n_1}$. Note
that $\tilde{f}_0^{n_0}$ and $\tilde{f}_1^{n_1}$ are periodic
functions. Arguing as in Theorem~\ref{thm:dec-real-line}, there is
a decomposition in compacts sets on the real line of
$L(\IFS(\tilde{f}_0^{n_0},\tilde{f}_1^{n_1}))
\setminus\{\pm\infty\}$. From the periodicity of
$\tilde{f}_i^{n_i}$, it follows that
each of these compact sets project on the circle as either, a periodic point of $f_i$, or as a $**$-intervals for $** \in \{ss, su, uu\}$.
By Theorem~\ref{t.duminy_lemma-per} these intervals are transitive sets for $\IFS(\Phi^n)$.

What is left is to study the limit set of a point whose
$\omega$-limit (or $\alpha$-limit) contains $\pm \infty$ on the
real line. This can only happen if there is a cycle for
$\IFS(\Phi^n)$ different from an \mbox{$ss$-interval.} Then
Theorem~\ref{thmA} implies that $\mathbb{S}^1$ is minimal for
$\IFS(\Phi^n)$ and $\IFS(\Phi^{-n})$, a contradiction. Therefore,
we obtain a decomposition of the limit set, as required.
\end{proof}

\section{A Denjoy type Theorem for IFS}

\subsection{A trichotomy for IFS.}

The next lemma shows some relations between the $\omega$-limit sets and the orbits of an IFS.
These properties will be necessary afterwards for the proof of the trichotomy result.
We will use the following notation: given a set $A$ on a manifold $M$,
we denote by $A'$ the set of accumulation points of $A$,
i.e., the set of points $y$ such that there exists a sequence
$(x_n)_n \subset A$ converging to $y$ with $x_n \not= y$ for all $n\in\mathbb{N}$.

\begin{lemma}
\label{esperanza}
Consider $\Phi=(\phi_1,\ldots,\phi_k) \in \mathrm{Hom}(M)^k$, a non-empty closed subset $K$ of a manifold $M$. Then it holds that:
\begin{enumerate}
\item \label{item32-1}
$\omega_\Phi(h(x))  \subset \omega_\Phi(x)\subset \overline{\mathrm{Orb}_\Phi^+(x)}$   for all  $x\in M$ and for $h\in \IFS(\Phi)$,
\item \label{item32-2}
if $K=\overline{\mathrm{Orb}^+_\Phi(x)}$ for all $x\in K$ then
$$
K= \phi_1(K)\cup\dots\cup\phi_k(K)=\omega_\Phi(x) \quad \text{for
all $x\in K$,}
$$
\item \label{item32-3}
$
\mathrm{Orb}_\Phi^+(x)\,'= \phi_1(\mathrm{Orb}_\Phi^+(x)\,'\,)\cup \dots \cup \phi_k(\mathrm{Orb}_\Phi^+(x)\,'\,)
$
for all $x\in M$.
\end{enumerate}
\end{lemma}
\begin{proof}By definition of  the $\omega$-limit set of a point
$x\in M$, $\omega_\Phi(h(x))\subset \omega_\Phi(x)$ for all $h\in
\IFS(\Phi)$. On the other hand, since
$$
 \overline{\mathrm{Orb}^+_\Phi(x)}=\{ y:  \text{there exists} \ (g_n)_n \subset \IFS(\Phi) \  \text{such that} \ y=\lim_{n\to \infty} g_n(x) \}
$$
then $\omega_\Phi(x)$ is a subset of the closure of the orbit of $x$. Therefore, we conclude~\eqref{item32-1}.

According to the first item, to obtain~\eqref{item32-2} it suffices to prove that
$$
 K \subset \bigcup_{i=1}^k \phi_i(K) \quad \text{and} \quad  K \subset \omega_\Phi(x) \  \ \text{for all
$x\in K$}.
$$
The first inclusion is followed from the minimality of $K$.
Indeed, for any $y\in K$ there is $(g_n)_n \subset \IFS(\Phi)$ such that $x=\lim_{n\to \infty} g_n(y)$.
Since we are working with a finitely generated semigroup, taking a subsequence if necessary,
we can assume that $g_n=\phi_i\circ \tilde{g}_n$ for some $i\in\{1,\ldots,k\}$ and $\tilde{g}_n \subset \IFS(\Phi)$.
Thus, $\phi^{-1}_i(x) \in K$ as it is the limit of $\tilde{g}_n(y)$ and $K$ is an invariant closed set.
Therefore $x=\phi_i(\phi^{-1}_i(x)) \in \phi_i(K)$, showing the first inclusion.

In order to prove the second inclusion, we fix $x,y\in K$ and consider a sequence of positive real numbers $\varepsilon_n=1/n \to 0$.
Let us construct by induction a sequence $(h_n)_n \subset \IFS(\Phi)$, such that the distance
between $y$ and
$h_n\circ \dots \circ h_1(x)$ is less than
$\varepsilon_n$. Since the orbit of $x$ is dense in $K$,
there is $h_1\in \IFS(\Phi)$ such that $d(y,h_1(x))<\varepsilon_1$.
Similarly, since $h_1(x)\in \mathrm{Orb}_\Phi^+(x) \subset K$ then  the orbit of $h_1(x)$ is dense in $K$
and there exists $h_2$ such that $d(y,h_2\circ h_1(x))<\varepsilon_2$.
Inductively we obtain the desired sequence $(h_n)_n\subset \IFS(\Phi)$. Hence,
$
y=\lim_{n\to \infty} h_n\circ\dots\circ h_1(x)
$
and thus $y\in \omega_\Phi(x)$ for all $x,y\in K$, concluding~\eqref{item32-2}.

Now we will show that $\mathrm{Orb}^+_\Phi(x)\,'$ is a self-similar set, that is we prove the last item.
Note that
$\phi_i(\mathrm{Orb}^+_\Phi(x)')\subset \mathrm{Orb}_\Phi^+(x)'$
for all $i=1,\ldots,k$. Indeed, if $y$ is an accumulation point of the orbit of $x$,
then $\phi_i(y)$ is approximated by points of the form $\phi_i\circ g_n(x) \in \mathrm{Orb}^+_\Phi(x)$
where $y=\lim_{n\to\infty} g_n(x)$ and $g_n(x) \not = y$ for all $n\in\mathbb{N}$.
This implies that $\phi_i(y)$ is also an accumulation point of $\mathrm{Orb}^+_\Phi(x)$ and so we conclude one of the inclusions.

To obtain the other inclusion
$$
   \mathrm{Orb}_\Phi^+(x)\,' \subset
\phi_1(\mathrm{Orb}^+_\Phi(x)\,'\,)
\cup\dots\cup \phi_k(\mathrm{Orb}^+_\Phi(x)\,'\,),
$$
we fix any $y\in \mathrm{Orb}_\Phi(x)\,'$. Then there exists $(g_n)\subset \IFS(\Phi)$
such that $y=\lim_{n\to\infty} g_n(x)$ and $g_n(x) \not = y$ for all $n\in\mathbb{N}$.
Since the semigroup $\IFS(\Phi)$ is finitely generated, taking a subsequence if necessary, we can assume that for some fixed
$i \in \{1,\ldots,k\}$ we have $g_n=\phi_i\circ \tilde{g}_n$ with $\tilde{g}_n \in \IFS(\Phi)$. Hence,
$
   \phi_i^{-1}(y)= \lim_{n\to \infty } \tilde{g}_n(x)
$ and $ \tilde{g}_n(y) \not = \phi^{-1}_i(y). $ Since the
accumulation set is a closed set,  it holds that
$\phi^{-1}_i(y)\in \mathrm{Orb}_\Phi^+(x)'$. This implies that  $y
\in \phi_i(\mathrm{Orb}^+_\Phi(x)')$ proving the desired inclusion
and therefore~\eqref{item32-3}. \end{proof}

The next trichotomy theorem is the IFS counterpart of the similar result for group actions
in \cite{Gh01}. An interval $I$ is forward-invariant for $\IFS(\Phi)$
if $f_i(I)\subset I$ for all $I$.

\begin{theo}{(Trichotomy for IFS)}
\label{p.IFS-tricotomi} Consider  $\IFS(\Phi)$  generated by
$\Phi=(f_1,\ldots,f_k)\in  \mathrm{Hom}(I)^k$, where $I$ is
$\mathbb{S}^1$ or a forward-invariant compact interval. Then there
exists a non-empty closed subset $K$ of $M$ such that
$$
K = \bigcup_{i=1}^k f_i(K) \quad \text{and} \quad \text{$K = \overline{\mathrm{Orb}^+_\Phi(x)}=w_\Phi(x)$ for all $x\in K$}.
$$
Moreover, one (and only one) of the following possibilities occurs:
\begin{enumerate}
\item $K$ is a finite orbit for $\IFS(\Phi)$,
\item $K$ has non-empty interior,
\item $K$ is a Cantor set.
\end{enumerate}
\end{theo}
\begin{proof}A family of non-empty closed subsets of $M$ such
that each member $\Lambda$ satisfies $ \Lambda=
f_1(\Lambda)\cup\dots\cup f_N(\Lambda) $ can be ordered by
inclusion. Since the intersection of nested compact sets is
compact and non-empty, Zorn's Lemma allows to conclude the
existence of a minimal (regarding the inclusion) non-empty closed
set $K$ such that $K=\cup f_i(K)$. We will show that
$K$ is an invariant minimal set for $\IFS(\Phi)$,
\begin{equation}
\label{eq:k-dens}
K=\overline{\mathrm{Orb}_\Phi^+(x)}\quad \text{for all $x\in K$.}
\end{equation}
Then, according to Item~\eqref{item32-2} in Lemma~\ref{esperanza}, we have $K=\overline{\mathrm{Orb}^+_\Phi(x)}=\omega_\Phi(x)$ for all $x\in K$.

In order to prove~\eqref{eq:k-dens},  since $K$ is an closed invariant set,
$
\mathrm{Orb}_\Phi^+(x)\,'\subset \overline{\mathrm{Orb}_\Phi^+(x)} \subset K
$
for all $x\in K$.
On other hand, according to~\eqref{item32-3} in Lemma~\ref{esperanza},  for each $x\in K$ the set of
accumulation points $\mathrm{Orb}^+_\Phi(x)\,'$ is a closed self-similar set.
Since  $K$ is minimal (regarding the inclusion) then  either, $\mathrm{Orb}^+_\Phi(x)\,'$ is an empty set or $K=\mathrm{Orb}_\Phi^+(x)\,'$.
We have two possibilities:  \begin{itemize}
\item[(i)] there is $x\in K$ such that $\mathrm{Orb}_\Phi^+(x)\,'$ is empty;
\item[(ii)] for all $x\in K$, it holds that $K=\mathrm{Orb}_\Phi^+(x)\,'$.
\end{itemize}
In the first case, it follows that $\mathrm{Orb}^+_\Phi(x)$ is a finite set, and therefore it is a non-empty closed self-similar set contained in $K$.
This implies, via Zorn's Lemma, that $K=\mathrm{Orb}_\Phi^+(x)=\omega_\Phi(x)$.

In the second case, we obtain that $K$ is an invariant minimal set for $\IFS(\Phi)$ and so
$K'=K$. Moreover, we have two options: $K$ has non-empty interior
or the interior of $K$ is empty and thus, since $M$ is a
one-dimensional manifold, $K$ is a Cantor set. This concludes the
proof of the theorem. \end{proof}

As an attractor for IFS is generated by a forward-invariant set, the theorem above
describes the shape of possible attractors of an IFS.

\begin{rem}\label{rem:Cantor}
In the case of $I$ is an $ss$-interval, which is forward-invariant, the end-points are global (with respect to $I$) attractors
and are necessarily part of any minimal set in $I$.
Thus the minimal set is unique. See also~\cite{S12}.
\end{rem}


\subsection{Proof of Theorem~\ref{thmC}.}

From Theorem~\ref{thmB} it follows that any minimal invariant closed set, $\Lambda$, for $\IFS(\Phi^n)$ is contained
in some piece of the spectral decomposition of the limit set. By the hypothesis of the theorem this piece
cannot be a single point because the generators do not have fixed points in common. It cannot also be
neither a $su$ nor a $uu$-interval because they are not forward invariant and any point can be made to leave the
interval. Thus necessarily the piece has to be a $ss$-interval $K^{ss}$, in which case by Remark~\ref{rem:Cantor} , the only  minimal invariant
set is $K^{ss}$ and so $\Lambda=K^{ss}$.

It is left to show that there cannot exist invariant Cantor sets
for $\IFS(\Phi)$. Let us suppose that $\Lambda$ is an minimal
invariant Cantor set for $\IFS(\Phi)$ and take any point
$p\in\Lambda$. Since $f_0$ and $f_1$ have no periodic points in
common, the basins of all attracting periodic points of $f_0$ and
$f_1$ cover $\mathbb{S}^1$. By hypothesis $\mathbb{S}^1$ is not
minimal and so from Theorem~\ref{thmA} there is at least one
$ss$-interval for $\IFS(\Phi^n)$ whose endpoints are attractors.
The cover of $\mathbb{S}^1$ by the basins of attractors allows $p$
to be attracted inside some $ss$-interval $K^{ss}$ for
$\IFS(\Phi^n)$. Since $\Lambda$ is a minimal invariant set for
$\IFS(\Phi)$, then $K^{ss}\subset \Lambda$, contradicting that it
has empty interior.

\subsection{A Theorem of Duminy on group actions of the circle.}

Finally, we will reprove Duminy's Theorem~\cite{Na11} under the assumption of Theorem~\ref{thmA}.
We denote by $\mathrm{G}(\Phi^n)$
the group generated by $\Phi^n=(f_0^{n_0},f_1^{n_1})$.

\begin{theo}[Duminy]
\label{t.Duminy}
There exists $\varepsilon>0{.}30$ such that if $f_0$ and $f_1$ are
diffeomorphisms on the circle with periodic points of period $n_0$ and $n_1$ respectively,
$\varepsilon$-close to the rotations in the $C^{1+\mathrm{bv}}$-topology and with no periodic points in common then,
\begin{itemize}
\item $\mathbb{S}^1$ is minimal for $\mathrm{G}(\Phi^n)$ and
\item the hyperbolic periodic points of  $\mathrm{G}(\Phi^n)$ are dense in $\mathbb{S}^1$.
\end{itemize}
\end{theo}
\begin{proof}Since there are no periodic points in common and we
are working with a group then considering the inverses if
necessary, there always exists a $ss$-interval and a cycle
different of the $ss$-interval for the group $\mathrm{G}(\Phi^n)$.
Moreover, this cycle has the same order as in
Lemma~\ref{claim:sk}. By Proposition~\ref{t.duminy_lemma-per}
every point in the $ss$-interval has a dense orbit for the action
of $\mathrm{G}(\Phi^n)$. Now similarly to the proof of minimality
in Theorem~\ref{thmB}, the cycle allows to bring any point and
interval to the $ss$-interval and this implies minimality of
$\mathbb{S}^1$ for $\mathrm{G}(\Phi^n)$.

The density of hyperbolic periodic points follows from the
argument described in Lemma~\ref{lem-densidad} (see also
Section~\ref{ss:density hyp-points}). This concludes the proof of
the result. \end{proof}

\section*{Acknowledgments}
Both authors thank  Enrique R. Pujals for useful
conversations and P. G. Barrientos is especially grateful to
Jos\'e A. Rodr\'iguez for the orientation and the helpful
comments.

During the preparation of this article the authors were partially
supported by the following fellowships: P. G. Barrientos by
FPU-grant AP2007- 031035, MTM2011-22956 project (Spain) and CNPq
post-doctoral fellowship (Brazil);  A. Raibekas by CNPq doctoral
fellowship (Brazil) and CAPES-PNPD (Brazil) pos-doctoral
fellowship. We thank Instituto Nacional de Matemática Pura e
Aplicada (IMPA) for its hospitality.

\bibliographystyle{plain}
\bibliography{br}

\end{document}